\newif\ifbka
\bkafalse 

\ifbka
\documentclass[lineno]{biometrika}
\RequirePackage{amssymb,mathtools}
\usepackage{nameref,zref-xr} 
\zxrsetup{toltxlabel} 
\zexternaldocument*{supplement}
\else
\documentclass{article}
\usepackage[margin=1.5in]{geometry}
\usepackage{ifthen}
\RequirePackage{amsmath,amsthm,amssymb,mathtools}

\newtheorem{proposition}{Proposition}
\newtheorem{lemma}[proposition]{Lemma}
\newtheorem{theorem}[proposition]{Theorem}
\newtheorem{assumption}{Assumption}
\newcommand{\tbl}{\caption}
\fi

\usepackage{bm}
\usepackage{natbib}

\RequirePackage{amsmath,amssymb,mathtools}
\usepackage{graphicx, multirow, placeins}
\usepackage{xr}

\DeclareMathOperator*{\argmin}{argmin}
\DeclareMathOperator*{\var}{var}

\DeclareMathOperator{\E}{E}

\DeclareMathOperator{\Pn}{P_n}
\let\P\relax 
\DeclareMathOperator{\P}{P}

\DeclarePairedDelimiter\abs{\lvert}{\rvert}
\DeclarePairedDelimiter\norm{\lVert}{\rVert}
\DeclarePairedDelimiter\set{\{}{\}}
\DeclarePairedDelimiter\inner{\langle}{\rangle}

\newcommand{\p}[1]{\left(#1\right)}
\newcommand{\sqb}[1]{\left[#1\right]}
\newcommand{\riesz}[1][]{\ifthenelse{\equal{#1}{}}{\gamma}{\gamma_{{#1}}}}
\newcommand{\hriesz}[1][]{\ifthenelse{\equal{#1}{}}{\hgamma}{\hgamma_{{#1}}}}
\newcommand{\ind}[1]{\mathbb I\p{#1}}
\newcommand{\bkaproof}{\ifbka%
\else %
Proof 
\fi
}

\newcommand{\mo}{m_0}
\newcommand{\hmo}{\hat{m}_0}

\newcommand{\he}{\hat{e}}

\newcommand{\hgamma}{\hat{\gamma}}
\newcommand{\hV}{\hat{V}}
\newcommand{\hL}{\hat{L}}

\newcommand{\F}{\mathcal{\model}}
\newcommand{\ff}{\mathcal{\model}}
\newcommand{\hh}{\mathcal{H}}

\newcommand{\xx}{\mathcal{X}}
\newcommand{\bb}{\mathcal{B}}

\newcommand{\treatment}{Z}
\newcommand{\target}{T}
\newcommand{\ntreatment}{n_{\treatment}}

\newcommand{\poptreatment}{\P_Z}

\newcommand{\influence}{\iota}
\newcommand{\model}{\mathcal{M}}

\newcommand{\HH}{\mathcal{H}}
\newcommand{\X}{\mathcal{X}}

\newcommand{\R}{\mathbb{R}}
\newcommand{\B}{\mathcal{B}}

\newcommand{\pr}{\mathrm{pr}}

\ifbka
\usepackage[plain,noend]{algorithm2e}

\makeatletter
\renewcommand{\algocf@captiontext}[2]{#1\algocf@typo. \AlCapFnt{}#2} 
\def\@algocf@capt@plain{top}
\renewcommand{\algocf@makecaption}[2]{%
  \addtolength{\hsize}{\algomargin}%
  \sbox\@tempboxa{\algocf@captiontext{#1}{#2}}%
  \ifdim\wd\@tempboxa >\hsize
    \hskip .5\algomargin%
    \parbox[t]{\hsize}{\algocf@captiontext{#1}{#2}}
  \else%
    \global\@minipagefalse%
    \hbox to\hsize{\box\@tempboxa}
  \fi%
  \addtolength{\hsize}{-\algomargin}%
}
\makeatother

\else

\title{Minimax Linear Estimation of the Retargeted Mean}
\author{David A. Hirshberg\thanks{\small Stanford University} \and Arian Maleki\thanks{\small Columbia University} \and Jos\'e R. Zubizarreta\thanks{\small Harvard University}}

\fi

\begin{document}

\ifbka
\jname{Biometrika}
\jyear{2020}
\jvol{103}
\jnum{1}
\accessdate{Advance Access publication on 31 July 2018}

\received{2 January 2017}
\revised{1 April 2017}

\markboth{Minimax Linear Estimation of the Retargeted Mean}{Minimax Linear Estimation of the Retargeted Mean}

\title{Minimax Linear Estimation of the Retargeted Mean}

\author{David A. Hirshberg}
\affil{Department of Statistics and Graduate School of Business, Stanford University,\\ Stanford, CA 94305, USA
\email{davidahirshberg@stanford.edu}}

\author{Arian Maleki}
\affil{Department of Statistics, Columbia University,\\ New York, NY 10027, USA
\email{arian@stat.columbia.edu}}

\author{\and Jos\'e R. ZUBIZARRETA}
\affil{Departments of Health Care Policy, Biostatistics, and Statistics, Harvard University, 180 Longwood Avenue, Boston, Massachusetts 02115, U.S.A. \email{zubizarreta@hcp.med.harvard.edu}}

\fi

\maketitle

\begin{abstract}
Evaluating treatments received by one population for application to a different target population of scientific interest is a central problem in causal inference from observational studies. We study the minimax linear estimator of the treatment-specific mean outcome on a target population and provide a theoretical basis for inference based on it. In particular, we provide a justification for the common practice of ignoring bias when building confidence intervals with these linear estimators.  Focusing on the case that the class of the unknown outcome function is the unit ball of a reproducing kernel Hilbert space, we show that the resulting linear estimator is asymptotically optimal under conditions only marginally stronger than those used with augmented estimators.  We establish bounds attesting to the estimator's good finite sample properties.  In an extensive simulation study, we observe promising performance of the estimator throughout a wide range of sample sizes, noise levels, and levels of overlap between the covariate distributions of the treated and target populations.
\end{abstract}

\ifbka
\fi

\section{Introduction}

The core challenge of causal inference is isolating the effect of a treatment from the confounding effect of differences between the characteristics of those who receive it and others. Often, this boils down to adjusting for the observed difference in covariate distributions between the target population for the application of a treatment and the group that we observe receiving it. Here we consider the estimation of the treatment-specific mean by such an adjustment, 
\[ \mu_w = \E\{m_w(X) \mid T=1\} \quad \text{ for } \quad m_w(x)=\E(Y \mid W=w, X=x), \]
where $W$ is a level of treatment, $X$ a vector of covariates, $Y$ the outcome of interest, and $T$ indicates membership in the target population.

For this task, minimax linear estimators and their variants, which use so-called \emph{balancing weights} to adjust for these baseline differences, have become increasingly popular \citep[e.g.,][]{armstrong2015optimal,fan2016improving,kallus2016generalized,wang2020minimal,wong2017kernel,zhao2019covariate,zubizarreta2015stable}. Given a model for the regression function, i.e., a set $\model$ containing $m_w$, the derivation of these estimators is straightforward: they are determined by a quadratic optimization problem.
Their use is supported by a long history of theoretical results \citep[e.g.,][]{donoho1994statistical,ibragimov1985nonparametric,juditsky2009nonparametric}. The usual approach to inference in these models involves \emph{bias-aware confidence intervals}, which are meant to accommodate bias arising from residual differences in covariate distributions. These intervals are defined in terms of the model $\model$, which must be specified a-priori; however, there are strong impossibility results about the use of an estimated
model in this context \citep{armstrong2015optimal}, and in practice, it is difficult to have intuition to select this model.
Given the practical importance of having valid confidence intervals,
in this paper we study the behavior of these bias-aware confidence intervals 
in the widely-used class of Reproducing Kernel Hilbert Space (RKHS) models \citep{kallus2016generalized, wong2017kernel, zhao2019covariate}. 

In this paper, we characterize the minimax linear estimator's bias adaptively,
investigating how our choice of model $\model$ interacts with the problem's intrinsic parameters,
the outcome regression function $m_w$ and a generalization $g_{\psi}$ of the inverse propensity score, 
to produce bias. We show that in correctly specified RKHS models, 
this bias is essentially always second order;
the minimax linear estimator is asymptotically
normal with optimal variance;
and bias-aware confidence intervals are asymptotically equivalent to the (bias-unaware)
asymptotically optimal interval.
From a practical standpoint, this provides a theoretical basis for the widespread practice of ignoring bias when constructing confidence intervals, as discussed by \citeauthor{kallus2016generalized} (\citeyear{kallus2016generalized}, Appendix B).

Furthermore, we show that the estimator is adaptive 
in its higher order behavior. Within families of RKHS models,
e.g., the classes $\model_k$ of functions with $k$ square integrable derivatives,
we benefit when our intrinsic parameters are smoother, e.g., $k'>k$ times differentiable,
without modeling this additional smoothness. That is, we show that it is unnecessary
to use the stronger model $\model_{k'}$ to get these benefits, as performance with the 
model $\model_k$ is comparable. Practically, this adaptivity property alleviates some 
of the burden of choosing the `right' model.

Our analysis is based largely on an equivalent dual characterization 
of the minimax linear estimator over an RKHS ball
as a plug-in estimator using a kernel ridge regression estimate of the regression function $m$. 
In this way, we build on the literature on random design ridge regression work in Hilbert spaces \citep[e.g.,][]{fischer2017sobolev, hsu2012random, smale2007learning}   
by focusing on bias in what is essentially a generalization task: the problem of predicting outcomes on target population different from that from which the training sample for the regression was drawn.

\section{Framework}

\subsection{The Estimand}

We consider an observational study in which for each unit $i$ we observe a covariate vector $X_i$, a categorical treatment status $W_i \in {0 \ldots C}$, and an outcome $Y_i=Y_i(W_i) \in \R$, and assume that  as a function of $(X_i,W_i)$, we can calculate indicators $T_i=T(X_i,W_i)$ that mark units as members of a \emph{target group} of scientific interest. 
Our goal is to estimate the average, over this target group, of the potential outcome $Y_i(0)$ that they would have been experienced had they received the treatment of interest $W_i=0$.
We assume that missingness arises from a strongly ignorable mechanism \citep{rosenbaum1983central}; see Appendix~\ref{sec:identification} for details.

\subsection{The Estimator}
\label{sec:understanding-minimax-linear}
For simplicity, we focus on estimating
$\pr(T_i=1)\mu_0$, which can be written as
\[ \psi(m) = \E\{T_i m_0(X) \} \quad \text{ for } \quad m_0(x)=\E(Y \mid W=0, X=x). \]
Our estimator $\hat{\psi}_{ML}$ is a minimax linear estimator of a sample-average version of this estimand, $\tilde{\psi}(m) = n^{-1} \sum_{i=1}^n T_i \mo(X_i)$, conditional on the study design $(X_i, W_i)_{i \le n}$. 
We choose the weights that result in the best estimate of $\tilde{\psi}(m)$ of the form $n^{-1}\sum_{i=1}^n \ind{W_i=0} \gamma_i Y_i$ in the worst case over regression functions $\mo$ in an absolutely convex class $\ff$ and over conditional variance functions $\var(Y_i \mid X_i=x, W_i=w)$  bounded by a constant $\sigma^2$. These weights $\hriesz \in \R^n$ solve the following convex optimization problem.

\begin{equation}
\label{eq:ml-weights}
\hriesz = \argmin_{\riesz} I_{\ff}^2(\riesz) + \frac{\sigma^2}{n^2}\norm{\riesz}^2,\  I_{\ff}(\riesz) = \sup_{f \in \ff}\frac{1}{n}\sum_{i=1}^n \set{\ind{W_i=0}\gamma_i - T_i}f(X_i).
\end{equation} 
Here $I_{\ff}$ measures how well the $\gamma$-weighted average of a function $f(x)$ observed on the treated sample 
matches the average over the target for all functions $f \in \ff$. 
As a result the minimax linear weights enforce the sample balance condition $n^{-1}\sum_{i=1}^n \ind{W_i=0} \gamma_i f(X_i) \approx n^{-1}\sum_{i=1}^n T_i f(X_i)$ uniformly over $\ff$. 
And when this class $\ff$ is chosen appropriately, the weights $\ind{W_i=0}\hgamma_i$ converge in empirical mean square to our functional's Riesz representer: the unique square integrable function $\riesz[\psi](X_i,W_i)$ that satisfies the corresponding population balance condition $\E\{\riesz[\psi](X_i,W_i)f(X_i,W_i)\} = \E\{ T_i f(X_i,0)\}$ for all square integrable functions $f(x,w)$
\citep{hirshberg2018augmented}. 
\begin{equation}
\label{eq:riesz-targeted-mean}
\riesz[\psi](x,w) = \ind{w=0}g_{\psi}(x) \ \text{ for } g_{\psi}(x) = \frac{\pr(T_i=1 \mid X_i=x)}{\pr(W_i=0 \mid X_i=x)}.
\end{equation}
We often call weights like $\riesz[\psi](X_i,W_i)$ inverse probability weights, as they invert the probabilistic mechanism that assigns units to our treatment and target groups to ensure unbiasedness
of the oracle estimator $\psi_{\star}=n^{-1}\sum_{i=1}^n \riesz[\psi](X_i,W_i)Y_i$.

As a first step toward understanding the behavior of our estimator, we decompose its error into a bias term and a noise term. 
We will consider estimation of the sample-average version of our estimand, $\tilde{\psi}(m) = n^{-1}\sum_{i=1}^n T_i m(X_i, 0)$,
as the behavior of the difference $\tilde{\psi}(m) - \psi(m)$ is out of our hands. We decompose it as follows.
\begin{equation}
 \label{eq:linear-error-decomp}
\begin{split}
&\hat{\psi}_{ML} - \tilde{\psi}(m)  
= \frac{1}{n}\sum_{i=1}^n \ind{W_i=0} \hgamma_i Y_i - T_i \mo(X_i) \\ 
&\quad = \frac{1}{n}\sum_{i=1}^n \set{ \ind{W_i=0} \hgamma_i - T_i } \mo(X_i) + \ind{W_i=0}\hgamma_i \varepsilon_i, \ \varepsilon_i = Y_i - m_0(X_i).
\end{split}
\end{equation}
It is clear from this expression that what we are minimizing in \eqref{eq:ml-weights} to define our weights is the mean squared error
conditional on $(X_i,W_i)_{i \le n}$, maximized over $\mo \in \F$.
 
At the heart of our argument will be a characterization of our estimator $\hat{\psi}_{ML}$ as the
average, over our target subsample, of a kernel-ridge-regression estimate $\hat{m}(\cdot)$ of the regression function $\mo$ based on the
subsample $\{i : W_i=0\}$ of units receiving the treatment of interest. Specifically,
\begin{lemma}
\label{lemma:ridge-equivalence}
If $\ff$ is the unit ball of an RKHS with norm $\norm{\cdot}$, then
\begin{align}
&\frac{1}{n}\sum_{i=1}^n \ind{W_i=0} \hgamma_i Y_i = \frac{1}{n}\sum_{i=1}^n T_i \hmo(X_i)\  \text{ where } \label{eq:ridge-equivalence} \\
&\hgamma = \argmin_{\gamma} I_{\ff}^2(\riesz) + \frac{\sigma^2}{n^2}\norm{\riesz}^2, \ I_{\ff}(\riesz) = \sup_{f \in \ff}\frac{1}{n}\sum_{i=1}^n \set{\ind{W_i=0}\gamma_i - T_i}f(X_i); \label{eq:ridge-equivalence-weighting} \\
&\hmo = \argmin_{m} \frac{1}{\ntreatment} \sum_{i:W_i=0} \{Y_i - m(X_i)\}^2 + \frac{\sigma^2}{\ntreatment} \norm{m}^2, \  \ntreatment=\abs*{\{i:W_i=0\}}. \label{eq:ridge-equivalence-regression}
\end{align}
\end{lemma}
We use this duality result to characterize the conditional bias of our estimator as the bias of this ridge regression estimator,
conditional on the design $(X_i,W_i)_{i \le n}$, averaged over the target subsample. As the weight $\sigma^2/n_Z$ of
the penalty term in our ridge regression is small, $\hat{m}$  will be fairly unbiased estimator of a regression function $m$ in our RKHS.
Furthermore, what bias there is will be attenuated by averaging over the target sample to a degree that depends on the smoothness of ratio of 
covariate densities within these two groups, $\pr(X_i=x \mid T_i=1) / \pr(X_i=x \mid W_i=0) \propto g_{\psi}(x)$.

Using this argument to characterize our bias term and a result from \citet[Theorem 2]{hirshberg2018augmented} 
to characterize our noise term, we obtain simple sufficient conditions for our estimator to be asymptotically efficient:
in essence, if we take a ball in a RKHS $\hh$ to be our model, it suffices that our regression function $\mo$ be in $\hh$.


\subsection{Reproducing Kernel Hilbert Spaces}
\label{sec:rkhs-review}

We study $\norm{\cdot}_{\hh}^2$-penalized regression, where $\hh$ is a RKHS of functions on a compact metric space $\xx$.
A RKHS is a Hilbert space on which point evaluation is continuous, i.e., there is a constant $C_x$ for each $x \in \xx$ such that $f(x) \le C_x \norm{f}_{\hh}$. For each $x \in \xx$, there is a \emph{representer} $K_x \in \hh$ satisfying $\inner{K_x,f}_{\hh}=f(x)$ for all $f \in \hh$. 
We call the function $K(x,y) = \inner{K_x, K_y}$ the \emph{kernel} of $\hh$.

Given any finite measure $\nu$ with support equal to $\xx$, we can completely characterize $\hh$ in terms of the operator
$L_{K,\nu} : L_2(\nu) \to L_2(\nu)$ defined as $(L_{K,\nu}f)(\cdot) = \int K(\cdot,x) f(x) d\nu(x)$.
By Mercer's theorem, its eigenvalues $\mu_j$ are nonnegative and summable, its eigenfunctions $\phi_j$ form an orthonormal basis for $L_2(\nu)$, and its scaled eigenfunctions
$\smash{\mu_j^{1/2}} \phi_j$ form an orthonormal basis for $\hh$
\citep[see, e.g.,][Chapter 4]{cucker2007learning}.
This establishes an isomorphism between $\hh$ and the space of square-summable sequences $\ell_2$. 
For every sequence \smash{$t \in \ell_2$}, we define the corresponding element in $\hh$ by
\smash{$f_t(x) = \sum_j t_j \mu_j^{1/2} \phi_j(x)$}, \smash{$\inner{f_t, f_u}_{\hh} = \inner{t,u}_{\ell_2}$.} We generalize this by
considering the family of spaces \smash{$[\hh]_{\nu}^{s}$}
induced by the sequence-space inner products \smash{$\inner{g_t, g_u}_{[\hh]_{\nu}^s} = \inner{t,u}_{\ell_2}$ for $g_t(x) = \sum_j t_j \mu_j^{s/2} \phi_j(x)$.}

The Sobolev spaces $H_k$ of $k$-times differentiable periodic functions 
on the unit cube in $\R^d$ endowed with Lebesgue measure $\mu$ are such a family --- the unit ball of $H_k$ is 
\smash{$\set{ \sum_{j \in \mathbb{Z}^d} t_j (1+\norm{j}_2^2)^{-k/2} e^{2\pi i j} : \sum_{j \in \mathbb{Z}^d} t_j^2 \le 1}$}
\citep[see][]{kuhn2014approximation}. $H_k$ is an RKHS for $k > d/2$ and for each such $k$ we get the same family of spaces:
$[H_k]_{\mu}^s = [H_{k'}]_{\mu}^{s'}$ if $ks=k's$.

\section{Main Results}
\label{sec:targeted-mean-main-results}
\subsection{Setting}
\label{sec:setting}
Throughout, we work in a setting in which we have independent and identically distributed observations $(X_i,W_i,Y_i)_{i \le n}$. For some binary function $T(x,w)$, we define an indicator for membership in the target population, $T_i=T(W_i,X_i)$.
Furthermore, we let $\poptreatment$ denote the distribution of the covariate vector $X_i$ on the subpopulation for which $W_i=0$,
and assume that $X_i$ is in a compact metric space $\xx$. 
As our model $\model$, we use the unit ball of a RKHS $\hh$ of functions on $\xx$ that is dense in $L_2(\poptreatment)$,
has a bounded kernel $K$, and for which \smash{$\norm{\cdot}_{[\hh]_{\poptreatment}^{\alpha}} \le A_{\alpha}\norm{\cdot}_{L_{\infty}(\poptreatment)}$}
for $\alpha < 1$ and some finite constant $A_{\alpha}$. 

\subsection{Asymptotic Results}
\begin{theorem}
\label{theo:asymptotic}
Suppose $\pr(W_i=0)>0$ and 
the target and treatment groups overlap in the sense that
the inverse propensity weight $g_{\psi}$ defined in \eqref{eq:riesz-targeted-mean} is bounded. Furthermore, suppose that $\mo \in \hh$ and either $\mo \in  [\hh]_{\poptreatment}^{1+\epsilon}$
or $g_{\psi} \in [\hh]_{\poptreatment}^{\epsilon}$ for $\epsilon > 0$.
Then the estimator 
\begin{equation} 
\label{eq:hpsi-ml}
\hat{\psi}_{ML} = n^{-1}\sum_{i=1}^n \ind{W_i=0} \hgamma_i Y_i 
\end{equation}
with weights $\hgamma$ defined in \eqref{eq:ml-weights} for any $\sigma > 0$ satisfies 
\begin{equation}
\label{eq:scaled-targeted-asymptotically-linear}
\begin{split}
& \hat{\psi}_{ML} - \psi(m) = n^{-1}\sum_{i=1}^n \influence(W_i,X_i,Y_i) + o_p(n^{-1/2}), \\
& \influence(w,x,y) = T(x,w)\mo(x) - \psi(m) + \riesz[\psi](x,w)\{y - \mo(x)\}.
\end{split}
\end{equation}
Furthermore, if $\pr(T_i=1) > 0$, then for $\hat p_T = n^{-1}\sum_{i=1}^n \ind{T_i=1}$,
\begin{equation}
\label{eq:targeted-asymptotically-linear}
\begin{split}
& \hat p_T^{-1} \hat{\psi}_{ML} - \mu_0 = n^{-1}\sum_{i=1}^n \influence'(W_i,X_i,Y_i) + o_p(n^{-1/2}), \\
& \influence'(w,x,y) = [T(w,x) (m_0(x) - \mu_0) + \riesz[\psi](x,w)\{y - \mo(x)\}] / P(T_i=1).
\end{split}
\end{equation}
\end{theorem} 
Thus, $n^{1/2}(\hat p_T^{-1}\hat{\psi}_{ML}- \mu_0)$ is asymptotically normal with mean zero and variance
$V =$\ \smash{$\E \{ \influence'(W_i,X_i,Z_i)^2$} \}, and given a consistent estimate \smash{$\hV$} of $V$,
$\hat p_T^{-1}\hat{\psi}_{ML} \pm z_{\alpha/2}\hV^{1/2}/n^{1/2}$
is an asymptotically valid confidence interval of level $1-\alpha$. 


This result offers a first-order characterization of the minimax linear estimator. It relies on the property that
the conditional bias term in our error decomposition \eqref{eq:linear-error-decomp} is second order.
In the next section, we will characterize this bias term more precisely and bound the correction used to form bias-aware intervals. This allows us to establish asymptotic equivalence between the bias-aware intervals and the interval suggested by Theorem~\ref{theo:asymptotic}.

\subsection{Finite Sample Results}
\label{sec:finite-sample-results}
Our finite sample results are stated in terms of the sequence of decreasing eigenvalues $\mu_1,\mu_2,\ldots$
of the kernel smoothing operator $L_{K,\poptreatment}$ defined in Section~\ref{sec:rkhs-review}.
We write $n_Z,p_Z$ for the number of units with $W_i=0$ and corresponding marginal probability; 
$n_T,p_T$ analogously for $T_i=1$; $n_\sigma = n p_Z / (2\sigma^2)$; and each instance of $c$ is a universal constant.

\begin{theorem}
\label{theo:imbalance}
Suppose that $g_{\psi} \in [\hh]_{\poptreatment}^{\kappa_g}$ for $\kappa_g \in [0,1]$. Then with probability $1-2q(\delta)$, 
\begin{equation}
\label{eq:uniform-imbalance-bound}
\begin{aligned}
I_{[\hh]_{\poptreatment}^{\kappa_m}}(\hgamma) &\le 
 c \delta^{-1} n^{-1/2} n_{\sigma}^{\frac{1-\kappa_m}{2}}\set*{\sum_{j=1}^{\infty} \min(\mu_j, n_{\sigma}^{-1})}^{1/2} 
+ 12 p_Z \norm{g_{\psi}}_{L_2(\poptreatment)} n_\sigma^{-\frac{\kappa_m + \kappa_g - \kappa_g(1-\kappa_m)}{2}} 
\end{aligned}
\end{equation}
for all $\kappa_m \in [1,2]$. Here,  letting $N_{\poptreatment}(x) = \sum_{j=1}^{\infty} \mu_j / (\mu_j + x)$, 
\[ q(\delta) = \delta + \exp\p{-\frac{np_Z}{10}} 
+ \exp\p{-\frac{(n p_Z)^{1-\alpha}2^{\alpha-4}\sigma^{2\alpha}}{\log[2e A_{\alpha}^2 \{1 + (\mu_1 n_\sigma)^{-1}\} N_{\poptreatment}(1 / 3 n_\sigma)]}}. \] 
\noindent
If $\mu_j \le cj^{-1/p}$ for $p < 1$ and $\mu_1 \ge c / n_{\sigma}$, it follows that with probability $1-2q(\delta)$, 
\begin{equation}
\label{eq:uniform-imbalance-bound-simplified}
\begin{aligned}
I_{[\hh]_{\poptreatment}^{\kappa_m}}(\hgamma) &\le
c \delta^{-1} (p n)^{-1/2} n_{\sigma}^{\frac{p-\kappa_m}{2}} + 12 p_Z \norm{g_{\psi}}_{L_2(\poptreatment)} n_\sigma^{-\frac{\kappa_m + \kappa_g - \kappa_g(1-\kappa_m)}{2}}
\end{aligned}
\end{equation}
for all $\kappa_m \in [1,2]$ where $q(\delta) \le  \delta + 2\exp\p{-c\sigma^{2\alpha}(n p_Z)^{1-\alpha}/\log[\max\{2, A_{\alpha},\ n_{\sigma}, 1/(1-p)\}]}$.
\end{theorem}

\subsection{Bias-aware confidence intervals}
\label{sec:bias-aware}
Bias-aware confidence intervals are formed by 
adding to the usual normal-approximation-based interval some additional slack to accomodate 
the estimator's maximal conditional bias
when $\norm{m_0}_{\hh}$ is bounded by some tuning parameter $\nu$:
\[ \hat \psi \pm \set{ z_{\alpha/2}n^{-1/2}\hat V^{1/2} + \nu I_{\hh}(\hgamma)}. \]
The slack factor $I_{\hh}(\hgamma)$ is bounded by the right side of \eqref{eq:uniform-imbalance-bound} with $\kappa_m=1$. 
This is negligible, so we get an interval first-order equivalent to
the usual one, if the inverse propensity score is at least infinitesimally smooth in the sense that $g_{\psi} \in [\hh]_{\poptreatment}^{\kappa_g}$ for $\kappa_g > 0$.
It is negligible even without this if we undersmooth, favoring control of bias over variance in our minimax approach by taking the noise level parameter $\sigma$ to zero as sample size increases.

Undersmoothing improves our control of bias, but it raises other concerns. It violates the stated assumptions of Theorem 1, which requires $\sigma$
to be nontrivially large to ensure 
convergence of the weights $\hriesz$ to the Riesz representer
$\riesz[\psi]$. However, this assumption can be relaxed. In Appendix~\ref{sec:asymptotics-proofs}, we show that for any RKHS $\hh$, this convergence happens for sequences $\sigma$ converging to zero slowly enough. Furthermore, taking $\sigma$ to zero too fast does not necessarily cause inferential problems.
The more we undersmooth, the closer a studentized version of $\hat\psi-\tilde\psi(m)$ 
will be to a studentized version of the noise term in our error decomposition \eqref{eq:linear-error-decomp}. Self-normalized Berry-Esseen bounds \citep[e.g.,][Theorem 5.9]{pena2008self}
imply this will be approximately normal irrespective of whether the weights $\hgamma$ converge, 
so long as individual weights are not too extreme.

\subsection{Adaptivity}
\label{sec:adaptivity}
The conditional bias bound \eqref{eq:uniform-imbalance-bound-simplified} depends both on a measure $p$ of the smoothness of the model space $\hh$ and measures $\kappa_{m}(\hh)$ and $\kappa_g(\hh)$ of the smoothness of $m_0$ and $g_{\psi}$ relative to $\hh$. It shows that our estimator is adaptive in the sense that it benefits from smoothness of $m_0$ and $g_{\psi}$ in excess of the minimal levels $\kappa_m = 1$ and $\kappa_g = 0$ implicit in our use of the model $\hh$ in our minimax framework. In fact, it suggests that substituting a stronger model \smash{$\hh_s = [\hh]_{\poptreatment}^s$} for $s>1$ can make things worse, even if it is correctly specified in the sense that $\mo \in \hh_s$.

Intuitively, increasing the smoothness of this space increases the strength of regularization; to compensate, we must decrease it another way, by taking $\sigma \to 0$. This is almost a complete wash with appropriate tuning of $\sigma$:
as $\smash{p(\hh_s)}=p(\hh)/s$ and $\smash{\kappa_f(\hh_s)} = \kappa_f(\hh)/s$ for $f \in \{m,g\}$,
this happens when $n_{\sigma} \propto \smash{n^s}$, i.e., for $\sigma^2 \propto \smash{n^{1-s}}$. 
The minimax linear estimator for $\hh_s$ and this smaller choice of $\sigma$ will generally be efficient when $m_0 \in \hh_s$, but our bound does not give a reason to prefer this approach to using a less-smooth model and $\sigma \propto 1$.
See Appendix~\ref{sec:asymptotics-proofs} for details.

\subsection{Invariance}
\label{sec:invariance}
The targeted treatment-specific mean $\mu_0$
is translation invariant in the sense that
if we increase $m_0(x)$ by a constant $t$,
$\mu_0$ increases by $t$. This means that 
the way outcomes are centered does not
impact comparisons between treatments.
The estimator \smash{$\hat p_T^{-1}\hat{\psi}_{ML}$} is not translation invariant, but modifying it to have this desirable property is straightforward. One such variant replaces $\hat{\psi}_{ML}$ with the estimator \smash{$\hat{\psi}_{ML_t}$} $= \hat p_T \bar{Y}_0 + n^{-1}\sum_{i:W_i=0} \hgamma_i (Y_i - \bar{Y}_0)$
where \smash{$\bar{Y}_0 = n_Z^{-1}\sum_{i:W_i=0} Y_i$.} This is a very simple augmented minimax linear estimator \citep{hirshberg2018augmented} 
incorporating a constant estimate \smash{$\bar{Y}_0$ of $\mo$.} See \citet[Section 4.5]{kallus2016generalized} for an alternative approach to invariance that substitutes an invariant seminorm for the norm $\norm{\cdot}_{\hh}$.



\section{Empirical Performance}
\label{sec:empirical-performance}

We evaluate our estimator on an example of \citet*{hainmueller}.
In this example, treatment is binary and the target is the whole population: $W_i \in \set{0,1}$ and $T_i=1$ for all $i$.
We observe $X_i \in \R^6$ with $X_{i1} \ldots X_{i3}$ normal with mean zero and covariance matrix $\Sigma$ defined in Appendix~\ref{sec:additional-simulations}, $X_{i4} \sim \operatorname{Uniform}([-3,3])$, $X_{i5} \sim \chi^2_1$, and $X_{i6} \sim \operatorname{Bernoulli}(1/2)$ independent of each other and of $X_{i1} \ldots X_{i3}$; missingness follows a probit model $pr(W_i = 0 \mid X_i) = \Psi\{\eta^{-1}(X_{i1} + 2 X_{i2} - 2 X_{i3} - X_{i4} - 0.5 X_{i5} + X_{i6})\}$; and outcomes follow a quadratic model $Y_i = (X_{i1} + X_{i2} + X_{i5})^2 + \sigma_{\varepsilon} \varepsilon_i$ for standard normal $\varepsilon_i$. 

The estimators compared are (i) an averaged regression estimator \smash{$n^{-1}\sum_{i=1}^n \hat{m}(X_i)$}, where \smash{$\hat{m}$} estimates $\mo(x)$ by ordinary least squares (OLS) on the treated units; (ii) an inverse propensity weighting (IPW) estimator \smash{$n^{-1}\sum_{i=1}^n \ind{W_i=0} Y_i / \he(X_i)$}, where \smash{$\he(x)$}
is a logistic regression estimator of $\pr(W_i = 0 \mid X_i=x)$; (iii) an augmented inverse probability weighting (AIPW) estimator 
\smash{$n^{-1}\sum_{i=1}^n \hat{m}(X_i) + \ind{W_i=0}\{Y_i - \hat{m}(X_i)\}/\he(X_i)$} incorporating the aforementioned estimators \smash{$\hat{m}$} and \smash{$\he$};
(iv) the translation invariant variant \smash{$\hat{\psi}_{ML_t}$} of the minimax linear estimator for the unit ball $\model$ of the RKHS with the Mat\'ern kernel
$K_{\nu}$ for $\nu=3/2$. We take $\sigma=0.1$. See Appendix~\ref{sec:additional-simulations} for details.


In addition to bias and root mean squared error, we will look at the width and coverage of 95\% confidence intervals of the form $\hat{\psi} \pm z_{.025} \hV^{1/2}/ n^{1/2}$,
where we use the 
variance estimator
$\hV = n^{-1}\sum_{i:T_i=1} \{\hat m_0 (X_i) - \hat{\psi}\}^2 + n^{-1}\sum_{i:W_i=0} \hriesz_i^2 \{Y_i - \hat m_0(X_i)\}^2$.
Here, $\hriesz_i$ are the weights used in the corresponding estimator and $\hat m_0$ is an OLS estimate of $\mo(x)$ based on the sample receiving treatment $W_i=0$.\footnotemark 
\footnotetext{The OLS estimator is not typically considered a weighting estimator, but it is linear in $Y$ and can therefore be expressed in that form.
Lemma~\ref{lemma:ridge-equivalence} shows that it is, in fact, a limiting ($\sigma \to 0$) case of our estimator $\hat{\psi}_{ML}$ in which we work with the RKHS of linear functions $f(x) = f^T x$
with the Euclidean inner product $\inner{f(x),g(x)}=f^T g$.}  

Table~\ref{table:sims} presents the results.
The logit missingness model used in the IPW and AIPW methods is barely misspecified, so the IPW and AIPW estimators perform well in moderate and large samples. Misspecification of the linear outcome model relied upon by OLS causes substantial bias in all sample sizes. 
The estimator $\hat{\psi}_{ML_t}$ performs well in all samples sizes,
and does so also in several variations on this example and in the classic example of \citet*{kang2007demystifying}.
We discuss these in Appendix~\ref{sec:additional-simulations}.

\begin{table}
\def~{\hphantom{0}}
\tbl{Performance in the example of Hainmueller with $\eta=\sqrt{30}$ and $\sigma_{\varepsilon}=10$.}{
\begin{tabular}{ccc}
Root-mean squared error & Bias & Coverage of 95\% interval \\
\begin{tabular}{l rrr}
n    & 50 & 200 & 1000  \\
IPW  & 4.37 & 3.86 & 1.21 \\
AIPW & 5.89 & 5.49 & 1.74 \\
OLS  & 4.78 & 3.19 & 2.66  \\ 
$\text{ML}_t$  & 3.30 & 1.68 & 0.94 
\end{tabular}
&
\begin{tabular}{l rrr}
n & 50 & 200 & 1000  \\
IPW  & -0.84 & -0.36 & -0.35 \\
AIPW & -1.09 & -0.48 & -0.48 \\
OLS  & -2.30 & -2.49 & -2.5  \\
$\text{ML}_t$  & 0.26 & -0.28 & -0.53  
\end{tabular}
&
\begin{tabular}{l rrr}
n & 50 & 200 & 1000 \\
IPW  & 0.93 & 0.95 & 0.99  \\ 
AIPW & 0.87 & 0.93 & 0.94  \\ 
OLS  & 0.93 & 0.87 & 0.46  \\
$\text{ML}_t$  & 0.93 & 0.97 & 0.99  
\end{tabular}
\end{tabular}
}
\label{table:sims}
\end{table}
\FloatBarrier

\ifbka
\bibliographystyle{biometrika}
\else
\bibliographystyle{plainnat-abbrev}
\fi
\bibliography{references}

\newpage

\ifbka
\else
\begin{appendix}
\section{Causal Identification}
\label{sec:identification}
The estimand we discuss in this paper arises frequently in causal inference. 
Consider an observational study in which for each unit $i$ we observe a covariate vector $X_i$, a categorical treatment status $W_i \in {0 \ldots C}$, and an outcome $Y_i=Y_i(W_i) \in \R$, and assume that  as a function of $(X_i,W_i)$, we can calculate indicators $T_i=T(X_i,W_i)$ that mark units as members of a \emph{target group} of scientific interest. 
Our goal is to estimate the average, over this target group, of the potential outcome $Y_i(0)$ that they would have been experienced had they received the treatment of interest $W_i=0$. 

The well-known problem of estimating a mean outcome when some outcomes are missing is such a problem.
In that case, we observe the outcome of interest if we observe an outcome at all and our target group is the entire population; i.e., we have $W_i=0$ if and only if we actually observe the outcome $Y_i$ and $T_i=1$ for all $i$.
However, the flexibility afforded us in this framework to define our target group introduces little additional complexity and can be valuable.
For example, if we are wondering whether to recommend a change to treatment $W_i=0$ for those who are above a given age and currently taking another treatment $W_i=1$, it is natural to estimate the average outcome we would expect to see for that specific group if that recommendation were followed, which we can do by defining $T_i$ in terms of both $X_i$ and $W_i$.

The mean with outcomes missing is identifiable when missingness arises from a strongly ignorable mechanism \citep{rosenbaum1983central}.
For our more general problem, we make the following assumptions generalizing those that comprise strong ignorability.
\begin{assumption}[Overlap]
\label{assumption:overlap}
The covariate distribution on the target population is dominated by that of the treatment $W_i=0$ population; i.e.,
\[ \pr(X_i \in \cdot \mid T_i=1) \ll \pr(X_i \in \cdot \mid W_i = 0). \] 
\end{assumption}
\begin{assumption}[Unconfoundedness]
\label{assumption:uncondoundedness}
Conditional on the covariates, the potential outcome mean is the same for units in the treatment and target groups; i.e.,
\[ \E\{Y_i(0) \mid X_i, W_i=0\} = \E\{Y_i(0) \mid X_i, T_i=1\}. \] 
\end{assumption}
Under these assumptions, our causal estimand $\E[Y_i(0) \mid T_i = 1]$ is identified as a linear functional of
the regression of the observed outcome on covariates and treatment, i.e., it may be written
\begin{equation*}
\label{eq:psi-targeted-mean}
\mu_0 = \E\{m_0(X_i) \mid T_i=1\} \ \text{ where } \ m_w(x) = \E(Y_i \mid X_i=x, W_i=w).
\end{equation*}

\section{A Sketch of the Proof of Theorem~\ref{theo:imbalance}}
\label{sec:proof-sketch}
The maximal conditional bias \smash{$I_{[\hh]_{\poptreatment}^{\kappa_m}}(\hgamma)$} is the supremum of $n^{-1}\sum_{i=1}^n \set{ \ind{W_i=0}\hgamma_i - T_i}\mo(X_i)$ over 
\smash{$\mo \in [\hh]_{\poptreatment}^{\kappa_m}$}. Via Lemma~\ref{lemma:ridge-equivalence}, for each individual $\mo$,
this is equal to the averaged error $n^{-1}\sum_{i=1}^n T_i (\hmo - \mo)(X_i)$ of a purely theoretical ridge regression estimator 
based on noiseless observations, $\hmo = \argmin_m n_Z^{-1}\sum_{i:W_i=0} \{\mo(X_i) - m(X_i)\}^2 + \lambda \norm{m}^2$ for $\lambda=\sigma^2/n_Z$.
Letting $b_{\mo}=\hmo - \mo$, we write this as the sum of a population average and its deviation from it, 
\begin{equation}
\label{eq:bias-two-term-decomp}
\begin{aligned}
  n^{-1}\sum_{i=1}^n T_i b_{\mo}(X_i) &= \int T(x,w)b_{\mo}(x)d\P \\ &+  \set*{ n^{-1}\sum_{i=1}^n T_i b_{\mo}(X_i) - \int T(x,w)b_{\mo}(x)d\P}. 
\end{aligned}
\end{equation}
To bound these terms, we rely on an explicit characterization of the bias function $b_{\mo}$ as a vector in $\hh$: 
$b_{\mo}=-\lambda(\hat L + \lambda I)^{-1}\mo$ where
$\hat L$ is a sample-based approximation to the RKHS smoothing operator $L_{K,\poptreatment}$. Using this,
we show that with high probability, for all $\mo \in [\hh]_{\poptreatment}^{\kappa_m}$, this bias function is in
a class $\bb$ of functions satisfying bounds on \smash{$\norm{b}_{[\hh]_{\poptreatment}^{k}}$} for $k \in {0,1,2}$.
We use standard integral operator techniques \cite[e.g.,][]{cucker2007learning, fischer2017sobolev}.

We conclude by bounding our population average and deviation terms uniformly over $b \in \bb$.
To bound the deviation term, we use empirical process techniques. 
To bound the population average, we perform a change of measure that reveals the benefits of the smoothness of $g_{\psi}$:
via \eqref{eq:riesz-targeted-mean}, the integral $\int f(x)T(x,w)d\P$ is equivalent to $\int f(x) p_Z g_{\psi}(x) d\poptreatment$,
i.e., $p_Z$ times the inner product $\inner{ f, g_{\psi}}_{L_2(\poptreatment)}$. This allows is to exploit the
the familial relationship between $L_2(\poptreatment)$ and $\hh$ that we discussed in Section~\ref{sec:rkhs-review}.
For any approximation $\tilde{g} \in \hh$ to $g_{\psi}$, via Cauchy-Schwarz,
\begin{align*}
 \inner{ g_{\psi}, b}_{L_2(\poptreatment)} 
&= \inner{ g_{\psi}-\tilde{g}, b}_{L_2(\poptreatment)} + \inner{ \tilde{g}, b}_{L_2(\poptreatment)} \\
&= \inner{ g_{\psi}-\tilde{g}, b}_{L_2(\poptreatment)} + \inner{ \tilde{g}, L_{K,\poptreatment} b}_{\hh}, \\
&\le \norm{g_{\psi}-\tilde{g}}_{L_2(\poptreatment)}\norm{b}_{L_2(\poptreatment)} + \norm{\tilde{g}}_{\hh}\norm{L_{K,\poptreatment} b}_{\hh}
\end{align*}
Insofar as $g_{\psi}$ is smooth, $\tilde g - g_{\psi}$ is small, so 
we care not about the magnitude of the bias $b$ itself but about the magnitude of the \emph{smoothed bias} $L_{K,\poptreatment}b$,
i.e., about $\norm{L_{K,\poptreatment} b}_{\hh}=\norm{b}_{[\hh]_{\poptreatment}^2}$.

The motivation behind this two-term bound based on the approximation $\tilde{g}$ merits some explanation. 
Suppose that instead of ridge regression, we estimated $\mo$ using least squares regression in the span of the first $K \le n$ eigenvectors $\phi_1 \ldots \phi_K$
of $L_{K,\poptreatment}$. Then 
the bias function $b$ would be orthogonal to any function $\tilde{g}$ in this span, and while 
$\inner{ g_{\psi}-\tilde{g}, b}_{L_2(\poptreatment)} = \inner{ g_{\psi}, b}_{L_2(\poptreatment)}$, 
we are better off applying Cauchy-Schwartz to the former, and when we do, it is best to take $\tilde{g}$ to be the best
approximation to $g_{\psi}$ in that span \citep[see, e.g.,][Proof of Lemma A5]{newey2018cross}.
In essence, insofar as the function $g_{\psi}$ is smooth,
the nonsmooth components of the function $b$ we are averaging against it are irrelevant.
When we do ridge regression, we do not get strict orthogonality because our estimator does
shrinkage along all eigenvectors of $\hL \approx L_{K,\poptreatment}$. 
However, much of this bias is along the higher order eigenvectors of $L_{K,\poptreatment}$
and therefore increasingly orthogonal to the functions in $[\hh]_{\poptreatment}^k$ as $k$ increases. 
By working with the smooth approximation $\tilde{g}$ to $g_{\psi}$, we exploit this.

\section{A Proof of Theorem~\ref{theo:imbalance}}
\label{sec:proving-finite-sample-bounds}
\subsection{Summary}
In this section, we will prove our finite sample bound on the maximal conditional bias \smash{$I_{[\hh]_{\poptreatment}^{\kappa_m}}(\hgamma)$}
from Theorem~\ref{theo:imbalance}. The lemmas we use are proven in Section~\ref{sec:minimax-linear-lemmas}.

The maximal conditional bias \smash{$I_{[\hh]_{\poptreatment}^{\kappa_m}}(\hgamma)$} is the supremum of $n^{-1}\sum_{i=1}^n \set{ \ind{W_i=0}\hgamma_i - T_i}\mo(X_i)$ over 
\smash{$\mo \in [\hh]_{\poptreatment}^{\kappa_m}$}. Via Lemma~\ref{lemma:ridge-equivalence}, for each individual $\mo$,
this is equal to the averaged error $n^{-1}\sum_{i=1}^n T_i (\hmo - \mo)(X_i)$ of a purely theoretical ridge regression estimator 
based on noiseless observations, $\hmo = \argmin_m n_Z^{-1}\sum_{i:W_i=0} \{\mo(X_i) - m(X_i)\}^2 + \lambda \norm{m}^2$ for $\lambda=\sigma^2/n_Z$.
It follows that the maximimal conditional bias is the supremum over \smash{$\mo \in [\hh]_{\poptreatment}^{\kappa_m}$}
of $n^{-1}\sum_{i=1}^n T_i b_{\mo}$ for $b_{\mo} = \hmo - \mo$. In the subsections folllowing, we will show that
$b_{\mo}$ is in some set $\bb_r$ with high probability for such $\mo$, bound this supremum uniformly over $b_{\mo} \in \bb_r$,
and then simplify our bound when the eigenvalues $\mu_j$ decay polynomially.

\subsection{Notation}
\label{sec:notation}
In terms of the rank-one operator $K_{x} K_{x}^T$ mapping $\HH \to \HH$ by $[K_{x} K_{x}^T]f= K_x \inner{K_x, f}_{\HH}$,
we define \smash{$\hL = \ntreatment^{-1}\sum_{i:W_i=0} K_{X_i} K_{X_i}^T$}, the empirical version 
of our integral operator $L = L_{K,\poptreatment}$. We define regularized variants $\hL_{\lambda} = \hL + \lambda I$ 
and $L_{\lambda} = L + \lambda I$. In place of \smash{$\norm{f}_{[\HH]_{\nu}^{s}}$}, we will often write equivalently  
\smash{$\norm{L^{-s/2} f}_{L_2(\nu)}$} or \smash{$\norm{L^{-(s+1)/2} f}_{\hh}$},
and we will let $\norm{L}_{\HH \to \HH}$ denote the operator norm $\sup_{\norm{x}_{\HH}\le 1}\norm{Lx}_{\HH}$.
We write $\poptreatment f$ meaning $\int f(x,w)d\poptreatment$, 
$n_Z = \sum_{i=1}^n \ind{W_i=0}$, $p_Z = \pr(W_i=0)$, $\lambda = \sigma^2/n_Z$, use $c$ to denote a universal constant 
that may differ from instance to instance, and write $a \ll b$ or $a=o(b)$ meaning $a/b \to 0$
and $a \sim b$ or $a \propto b$ meaning there exist constants $L,U$ for which $L \le a/b \le U$.

For a set of functions $f(z)$
and $\sigma_1 \ldots \sigma_n$ independent with $\pr(\sigma_i=1)=\pr(\sigma_i=-1)=1/2$, we call 
$R_n(\F)=\E\{ \sup_{f \in \F} \abs{n^{-1}\sum_{i=1}^n \sigma_i f(Z_i)}\}$ 
the \emph{Rademacher complexity} of $\F$. We write $\F_r = \{ f \in \F : \E\{ f(Z_i)^2 \} \le r^2 \}$,
and call $R_n(\F_r)$ the \emph{local Rademacher complexity.}

\subsection{Characterizing $b_{\mo}$}

This optimization problem defining $\hmo$ has an explicit solution \citep[see e.g.][]{hsu2012random},
\[ \hmo= \hL_{\lambda}^{-1} \ntreatment^{-1}\sum_{i:W_i=0} K_{X_i} \mo(X_i) 
       = \hL_{\lambda}^{-1} \ntreatment^{-1}\sum_{i:W_i=0} K_{X_i} \inner{K_{X_i}, \mo}_{\hh} 
       = \hL_{\lambda}^{-1} \hL \mo. \]
Here $\hmo$ is written abstractly as a random element of the Hilbert space $\hh$: 
its evaluation at a point $x$ is $\hmo(x) = \inner{K_x, \hmo}_{\hh}$.
Thus, we may write our bias function, $b_{\mo}(x)=\hmo(x) - \mo(x)$, as
$\inner{K_x, b_{\mo}}$ where $b_{\mo} = [\hL_{\lambda}^{-1}\hL - I]\mo = -\lambda \hL_{\lambda}^{-1}\mo$.
The first step of our proof is to show that for any $\kappa_m \in [1,2]$, 
on an event of high probability, for all \smash{$\mo \in [\hh]_{\poptreatment}^{\kappa_m}$}, 
$b_{\mo}$ is in a set $\B_{r} = \set{ h \in \hh : \norm{L^{k/2}h}_{\hh} \le  r_k \ \text{ for }\ k \in \set{0,1,2}}$
with $r$ depending on $\kappa_m$ and $\kappa_g$.

To do this, observe that for any $\eta$,
\[ L^{\eta} b_{\mo} = -\lambda L^{\eta} \hL_{\lambda}^{-1} m_0 
                = -\lambda (L^{\eta}L_{\lambda}^{-\eta}) (L_{\lambda}^{\eta}\hL_{\lambda}^{-1} L_{\lambda}^{1-\eta}) 
		    (L_{\lambda}^{\eta-1}L^{\epsilon}) (L^{-\epsilon} \mo) \]
and consequently
\[ \norm{L^{\eta} b_{\mo}}_{\hh} \le \lambda \norm{L^{\eta}L_{\lambda}^{-\eta}}_{\hh \to \hh} 
					 \norm{L_{\lambda}^{\eta}\hL_{\lambda}^{-1} L_{\lambda}^{1-\eta}}_{\hh \to \hh}
					 \norm{L_{\lambda}^{\eta-1}L^{\epsilon}}_{\hh \to \hh} \norm{L^{-\epsilon} \mo}_{\hh}. \]
In this bound, the first and third operator norm factors are bounded by $1$ and $\lambda^{\epsilon + \eta - 1}$ respectively 
when $\epsilon +\eta \le 1$. To see the latter bound holds, observe that
$\norm{L_{\lambda}^{\eta-1}L^{\epsilon}}_{\hh \to \hh}=\norm{L_{\lambda}^{-1}L^{\epsilon/(1-\eta)}}^{1-\eta}_{\hh \to \hh}$,
take $\nu=\epsilon/(1-\eta)$ in following lemma, and observe that $\nu-1=(\epsilon+\eta - 1)/(1-\eta)$
and therefore $(\lambda^{\nu-1})^{1-\eta} = \lambda^{ \epsilon+\eta-1}$.
As $q(x) \ge 1$ on $[0,1]$, we use the simplified bound $\lambda^{\nu-1}$.

\begin{lemma}
\label{lemma:regularized-product-opnorm-bound}
Let $L$ be a compact operator on an RKHS $\hh$. Then if $\nu \in [0,1]$,
\[ \norm{L_{\lambda}^{-1} L^{\nu}}_{\hh \to \hh} \le  q(\nu)^{-1} \lambda^{\nu - 1}  \ \text{ where }
 \ q(x) = \p{\frac{x}{1-x}}^{1-x} + \p{\frac{x}{1-x}}^{-x}.\]
\end{lemma}

\citet[Equation 26]{fischer2017sobolev} bounds
 the remaining factor.
\begin{align*}
&\norm{L_{\lambda}^{\eta} \hL_{\lambda}^{-1} L_{\lambda}^{1-\eta}}_{\hh\to\hh} = \norm{L_{\lambda}^{1/2} \hL_{\lambda}^{-1} L_{\lambda}^{1/2}}_{\hh\to\hh} \le 3
\quad \text{ with probability}\ 1-2\exp(-u) \text{ if } \\
&n_Z\lambda^{\alpha} \ge 8u A_{\alpha}^2  \log\{2eN_{\poptreatment}(\lambda)(1+\lambda/\mu_1)\},\ 
N_{\poptreatment}(\lambda) = \sum_{j=1}^{\infty} \mu_j / (\mu_j + \lambda).
\end{align*}
When this holds, $\norm{L^{\eta} b_{\mo}}_{\hh} \le 3 \norm{L^{-\epsilon} \mo}_{\hh} \lambda^{\epsilon + \eta}$ for $\epsilon,\eta$ satisfying $\epsilon+\eta \le 1$.

To eliminate this bound's dependence on $\ntreatment$, observe that $\ntreatment$ is the sum of $n$ independent and identically distributed Bernoullis $\ind{W_i=0}$,
so by the multiplicative Chernoff bound it 
is in the range $[(1/2)np_Z, (3/2)np_Z]$ with probabiity $1-2\exp(-np_Z/10)$ \citep[see e.g.][Theorem 4.5]{mitzenmacher2005probability}.
And it follows that $\lambda = \sigma^2/n_Z$ is in the range $[\lambda_{\star}/3, \lambda_{\star}]$
for $\lambda_\star = 2\sigma^2/(np_Z)$. Because our bound on $\norm{L^{\eta} b_{\mo}}_{\hh}$ is increasing in $\lambda$,
on an event on which this Chernoff bound holds we get a valid bound by substituting $\lambda_\star$.
Furthermore, as $n_Z \lambda^{\alpha} = \sigma^{2\alpha} n_Z^{1-\alpha}$, on this event our lower bound on $n_Z \lambda^{\alpha}$ 
is satisfied if $\sigma^{2\alpha}(np_Z/2)^{1-\alpha} \ge 8u A_{\alpha}^2 \log\{2eN_{\poptreatment}(\lambda_\star/3)$ $(1+\lambda_\star/\mu_1)\}$.
Thus, with probability $1 -  2\exp(-np_Z/10) - 2\exp(-u)$, all \smash{$\mo \in [\hh]_{\poptreatment}^{\kappa_m}$} satisfy
\begin{equation}
\label{eq:b-containment}
\begin{aligned}
&\norm{L^{\eta} b_{\mo}}_{\hh} \le 3 \norm{L^{-\epsilon} \mo}_{\hh} \lambda_\star^{\epsilon + \eta} \ \text{ for }\ \epsilon+\eta \le 1 \ \text{ and }\ \\
&u = (np_Z)^{1-\alpha} 2^{\alpha - 4} \sigma^{2\alpha} A_{\alpha}^{-2} / \log\set{2eN_{\poptreatment}(\lambda_\star/3)(1+\lambda_\star/\mu_1)}.
\end{aligned}
\end{equation}
When this holds, $b_{\mo} \in \B_r$ with $r_k= 3 \norm{L^{-\epsilon} \mo}_{\hh}\lambda_\star^{\min(\epsilon + k/2,1)}$ for all $\mo \in \hh$.
Reparameterizing,
\[ \norm{L^{-\epsilon} \mo}_{\hh}=\norm{L^{-(1+2\epsilon)/2} \mo}_{L_2(\poptreatment)} = \norm{\mo}_{[\hh]_{\poptreatment}^{\kappa_m}}
 \ \text{ for } \kappa_m = 1+2\epsilon, \]
and $b_{\mo} \in \B_r$ with $r_k= 3 \lambda_\star^{\min\set{(\kappa_m + k - 1)/2,1}}$ for all $\mo \in [\hh]_{\poptreatment}^{\kappa_m}$
when $\kappa_m \ge 1$.

\subsection{Bounding conditional bias uniformly}
We bound $n^{-1}\sum_{i=1}^n T_i b$ uniformly over $b \in \B_r$ in two pieces. By the triangle inequality,
\[ \sup_{b \in \B_{r}} \abs*{n^{-1}\sum_{i=1}^n T_i b(X_i)} \le
    \sup_{b \in \B_{r}} \abs*{n^{-1}\sum_{i=1}^n T_i b(X_i) - \E T_i b(X_i)} + 
  \sup_{b \in \B_{r}} \abs*{\E T_i b(X_i)}. \]
The first term is the maximum of a mean-zero empirical process, 
and via the symmetrization and contraction inequalities \citep[e.g.,][Theorems 3.1.21 and 3.1.16]{gine2015mathematical}, its mean 
is bounded by $2R_n(\B_{r})$.
By Markov's inequality, this mean bound implies the probability $1-2\delta$ bound $\delta^{-1}R_n(\B_{r})$,
where $R_n(\B_{r}) \le cn^{-1/2}r_0[\sum_{j=1}^{\infty} \min\{\mu_j, (r_1/r_0)^2)\}]^{1/2}$ \citep[Theorem 41]{mendelson2002geometric}.
To bound the second term, observe that for any function $b$,
\begin{equation*}
\label{eq:change-of-measure-bound}
\begin{aligned}
\E\{T b(X)\} &= \E\{\ind{W=0} g_{\psi}(X) b(X)\} = p_Z\E\{g_{\psi}(X)b(X) \mid W=0\} \\
	   &= p_Z\p{ \inner{g_{\psi} - \tilde g, b}_{L_2(\poptreatment)} 
	    +		 \inner{\tilde g, b}_{L_2(\poptreatment)}} \\
	   &\le p_Z\p{ \norm{g_{\psi} - \tilde g}_{L_2(\poptreatment)}\norm{b}_{L_2(\poptreatment)} +
			  \norm{\tilde g}_{\hh}\norm{Lb}_{\hh}}.
\end{aligned}
\end{equation*}
Here we've used the Riesz representation \eqref{eq:riesz-targeted-mean} of the functional 
mapping $b \to \E(T b(X))$; expanded $g_{\psi}$ around a smooth approximation $\tilde g$;
and taken the Cauchy-Schwarz bound on each resulting term, using in the second case the identity
$\inner{g, Lb}_{\hh} = \inner{g, b}_{L_2(\poptreatment)}$.
This implies that 
$\E\{T b(X)\} \le  p_Z r_1 [ \norm{g_{\psi} - \tilde g}_{L_2(\poptreatment)} + (r_2/r_1)\norm{\tilde g}_{\hh} ]$ for all $b \in \B_{r}$.
As a result of Lemma~\ref{lemma:rkhs-approx-extremal} below, 
we can choose $\tilde g$ such that this implies $\E\{T b(X)\} \le  4 p_Z  r_1 (r_2/r_1)^{\kappa_g} \norm{g_{\psi}}_{L_2(\poptreatment)}$.
Adding our two bounds together, with probability $1-2\delta$,
\begin{equation}
\label{eq:uniform-bias-bound}
\begin{aligned}
\sup_{b \in \B_{r}} \abs*{n^{-1}\sum_{i=1}^n T_i b(X_i)} 
&\le c \delta^{-1} n^{-1/2} r_0 \sqb{\sum_{j=1}^{\infty} \min\{\mu_j, (r_1/r_0)^2)\}}^{1/2} \\
&+  4 r_1 p_Z \norm{g_{\psi}}_{L_2(\poptreatment)}  (r_2/r_1)^{\kappa_g}.
\end{aligned}
\end{equation}

\begin{lemma}
\label{lemma:rkhs-approx-extremal}
If $\norm{L_{K,\nu}^{\eta} g}_{L_2(\nu)} < \infty$ for $\eta \in [0,1/2]$, there exists $\tilde g$ such that 
\[ \norm{\tilde g - g}_{L_2(\nu)} + t \norm{\tilde g}_{\hh} \le t^{2\eta} \cdot 2 \norm{g}_{L_2(\nu)} q(2\eta) \ \text{ where } 
 q(x) = \p{\frac{x}{1-x}}^{1-x} + \p{\frac{x}{1-x}}^{-x}.\]
\end{lemma}
As $q(x) \le 2$ on $[0,1]$, we use the simplified bound $4t^{2\eta}\norm{g}_{L_2(\nu)}$. 

With probability $1 - 2\delta - 2\exp(-np_Z/10) - 2\exp(-u)$, both \eqref{eq:b-containment} and \eqref{eq:uniform-bias-bound} hold, and
\begin{equation}
\label{eq:bias-bound-general}
\begin{aligned}
\sup_{\mo \in [\hh]_{\poptreatment}^{\kappa_m}} \abs{n^{-1}\sum_{i=1}^n T_i b(X_i)} \le 
& c\delta^{-1} n^{-1/2}\lambda_\star^{\frac{\kappa_m - 1}{2}} \set*{\sum_{j=1}^{\infty} \min(\mu_j, \lambda_\star)}^{1/2} \\ 
+\quad &  12 p_Z \norm{g_{\psi}}_{L_2(\poptreatment)} \lambda_\star^{\kappa_m/2 + \kappa_g(1-\kappa_m/2)}.
\end{aligned}
\end{equation}
Here we've used the implication of \eqref{eq:b-containment} that for $\mo \in [\hh]_{\poptreatment}^{\kappa_m}$
for $\kappa_m \in [1,2]$, $b_{\mo} \in \B_r$ with 
\[ r_0 = 3 \lambda_\star^{(\kappa_m - 1)/2},\ r_1 = 3 \lambda_\star^{\kappa_m/2},\ r_2 = 3 \lambda_\star \ \text{ satisfying }\  
r_1/r_0 = \lambda_{\star}^{1/2}, \ r_2/r_1 = \lambda_{\star}^{1-\kappa_m/2}. \]
This is equivalent to \eqref{eq:uniform-imbalance-bound} from Theorem~\ref{theo:imbalance}, in which we rename $\lambda_\star$ to $n_{\sigma}^{-1}$
and rearrange the sum $\kappa_m/2 + \kappa_g(1-\kappa_m/2)=\{\kappa_m + \kappa_g + \kappa_g(1-\kappa_m)\}/2$.

When $\mu_j \le c j^{-1/p}$ for all $j$ with $p < 1$, we can derive from this a relatively simple bound.  
Given this bound on the eigenvalue sequence, as $\lambda \le x^{-1/p}$ if and only if $x \le \lambda^{-p}$, 
\begin{equation}
\label{eq:integral-approx}
\begin{aligned}
\sum_{j=1}^{\infty} \min(\mu_j, \lambda) 
&\le c\int_{0}^{\lambda^{-p}} \lambda + c\int_{\lambda^{-p}}^{\infty} x^{-1/p} \\
&= c\lambda^{1-p} + \{ c/(1-1/p) \} (\lambda^{-p})^{1-1/p} \le \{ c / (1-1/p) \} \lambda^{1-p}.
\end{aligned}
\end{equation}
Furthermore, $N_{\poptreatment}(\lambda) \le (c/\lambda)^p/(1-p)$ \citep[comments below Equation 14]{fischer2017sobolev}. 
Making these substitutions yields \eqref{eq:uniform-imbalance-bound-simplified}.

\section{Proof of Asymptotic Results}
\label{sec:asymptotics-proofs}

\subsection{Summary}
\label{sec:asymptotics-summary}
In this section, we prove a generalization of Theorem~\ref{theo:asymptotic} that allows $\sigma$ to vary with $n$.
We generalize Theorem~\ref{theo:asymptotic} by replacing its second sentence with the following assumption.
Througout, we use the notation described in Section~\ref{sec:notation}.

\begin{assumption}
\label{assu:sigma}
For some $\kappa_m \in [1,2]$ and $\kappa_g \in [0,1]$,
\smash{$\norm{m_0}_{[\hh]_{\poptreatment}^{\kappa_m}}$} and \smash{$\norm{g_{\psi}}_{[\hh]_{\poptreatment}^{\kappa_g}}$} are finite
and $\sigma$ satisfies, for $\kappa = \kappa_m + \kappa_g + \kappa_g (1-\kappa_m)$ and any $\epsilon > 0$,
\begin{enumerate}
\item[(i)]  $n^{1-1/\kappa} \gg \sigma^2 \ge n^{1+\epsilon- 1/\max(\alpha, p)}$ if $\mu_j \le cj^{-1/p}$ for $p<1$; 
\item[(ii)] $n^{1-1/\kappa} \gg \sigma^2 \gg a_n$ otherwise for some sequence $a_n \to 0$ depending on $\hh$.
\end{enumerate}
\end{assumption}
This generalizes Theorem~\ref{theo:asymptotic}, as under those assumptions $\kappa > 1$ and therefore constant $\sigma$
satisfies \smash{$n^{1-1/\kappa} \gg \sigma^2 \gg a_n$}. The additional flexibility to take $\sigma \to 0$ justifies the
claim from Section~\ref{sec:bias-aware} that we can get results analogous to Theorem~\ref{theo:asymptotic} while undersmoothing. 
Furthermore, it justifies the claim from Section~\ref{sec:adaptivity} that the estimator will be efficient for $\sigma^2 \propto n^{1-s}$ 
when we use the model \smash{$\hh_s$} for $s > 1$. To be precise about the latter,
let $\hh_s =[\hh']_{\poptreatment}^s$
where $\hh_1$ satisfies the assumptions in Section~\ref{sec:setting} and the associated eigenvalue sequence
satisfies the bound \smash{$\mu_j(\hh_1) \le c j^{-1/p_1}$} for $p_1<1$. Assumption~\ref{assu:sigma}
is satisfied for $\hh=\hh_s$ and $\sigma^2 = c n^{1-s}$
 when $s \ge 1$, $m_0 \in \hh_s$, and either $m_0 \in \hh_{s+\epsilon}$ or $g_{\psi} \in \hh_{\epsilon}$ for $\epsilon > 0$.
To see this, observe that these conditions on $m_0$ and $g_{\psi}$ imply the assumption's conditions involving $\kappa_m/\kappa_g/\kappa$, 
and the values of $\alpha_s$ and $p_s$ associated with $\hh_s$ satisfy $\alpha_s = \alpha_1/s$ and $p_s=p_1/s$.
Thus, the lower bound $\sigma^2=cn^{1-s} \ge n^{1+\epsilon- 1/\max(\alpha, p)}$
is equivalent to $s < 1/\max(\alpha,p)=s/\max(\alpha_1,p_1)$ and therefore to our assumption $\max(\alpha_1,p_1) < 1$.

In the subsections following, we will show that the bias term in our error decomposition \eqref{eq:linear-error-decomp}
is asymptotically negligible using Theorem~\ref{theo:imbalance}, show that the weights $\hriesz$ converge to the Riesz representer $\riesz[\psi]$
by a reduction to Theorem 2 of \citet{hirshberg2018augmented}, and use these results to show 
that $\hat\psi_{ML}$ and $p_T^{-1}\hat\psi_{ML}$ are asymptotically linear estimators of $\psi(m)=\pr(T_i=0)\mu_0$ and $\mu_0$ respectively
as claimed. As $\hh$ is fixed in our asymptotic regime, we will use $c$ to denote a constants that may depend on $\hh$
and which, as before, may vary from instance to instance.

We will use the following lemma.
\begin{lemma}
\label{lemma:complexity}
Let $\hh$ be an RKHS of functions on a compact metric space $\X$ associated with a kernel $K$
and unit ball $\B$, 
$\nu$ a measure with support equal to $\X$, 
$L$ the integral operator defined by $(L f)(\cdot) = \int K(\cdot, x)f(x)d\nu(x)$,
and $\mu_j$ the sequence of decreasing eigenvalues of $L$.
Then if $r_n \to 0$, the local Rademacher complexity 
$R_n(\B_{r_n})$ is $o(n^{-1/2})$ and 
the related fixed point $r_{\star} = \inf\{r > 0 : R_n(\B_{r}) \le cr^2\}$
is $o(n^{-1/4})$. Furthermore, if $\mu_j \le c j^{-1/p}$ for $p<1$, 
then $R_n(\B_{r}) \le cn^{-1/2}r^{2(1-p)}$ and $r_\star \le c n^{-1/(4p)}$.
\end{lemma}

\begin{proof}
We begin with the bound $R_n(\hh_{r}) \le c n^{-1/2} \{\sum_{j=1}^{\infty} \min(\mu_j, r^2)\}^{1/2}$
of \citet[Theorem 41]{mendelson2002geometric}. If $\mu_j \le cj^{-1/p}$ for $p < 1$, $R_n(\hh_{r}) \le cn^{-1/2}r^{2(1-p)}$
via the integral approximation \eqref{eq:integral-approx} and \smash{$r_\star \le c n^{-1/(4p)}$}, as \smash{$cn^{-1/2}r^{2(1-p)} \le cr^2$} 
when $r$ exceeds this bound. 

More generally, the aforementioned bound implies that for any integer $J$, $R_n(\hh_{r}) \le c n^{-1/2} \{J r^2 + \sum_{j > J} \mu_j\}^{1/2}$.
Because the eigenvalues $\mu_j$ are summable (see Section~\ref{sec:rkhs-review}), the tail sums $\sum_{j > J}\mu_j$ converge to zero as $J \to \infty$.
And if $r_n \to 0$ as $n \to \infty$, we can take $J_n \to \infty$ slow enough that $J_n r_n^2 \to 0$. It follows that
$R_n(\hh_{r_n}) = o(n^{-1/2})$ for $r_n \to 0$, as for such $J_n$ the bound above is $o(n^{-1/2})$.
Furthermore, $r_\star \le r$ if $r$ satisfies $Jr^2 \le cnr^4$ and $\sum_{j > J} \mu_j \le cnr^4$ for some $J$. 
This holds for $r_n=n^{-1/4}a_n$ for some sequence $a_n \to 0$, as for this choice of $r_n$, these conditions reduce to $J_n \le cn^{1/2}a_n^2$ and
$\sum_{j > J_n} \mu_j \le ca_n^4$.  To see this, observe that because tail sums converge to zero, 
for any sequence $J_n \to \infty$ satisfying $J_n \ll n^{1/2}$, 
there exists a sequence $a_n$ converging to zero slowly enough that $\sum_{j > J_n} \mu_j \le ca_n^4$
and $J_n \le cn^{1/2}a_n^2$. Thus, $r_{\star} = o(n^{-1/4})$.
\end{proof}

\subsection{Negligible bias}
\label{sec:negligble-bias}

Because $\mo \in \norm{\mo}_{\hh}\hh$, 
the conditional bias term in our error decomposition \eqref{eq:linear-error-decomp} 
is bounded by the maximal conditional bias over $\norm{\mo}_{\hh}\hh$,
\smash{$I_{\norm{\mo}_{\hh}\hh}(\hgamma) = \norm{\mo}_{\hh}I_{\hh}(\hgamma)$}. And because
we have assumed $\norm{\mo}_{\hh} < \infty$, it suffices to show that $I_{\hh}(\hgamma)=o_p(n^{-1/2})$.
We use \eqref{eq:uniform-imbalance-bound} from Theorem~\ref{theo:imbalance}.

First, we will show that $q(\delta) \to \delta$ as $n \to \infty$, so in large enough samples \eqref{eq:uniform-imbalance-bound}
holds with arbitrarily high probability. As $p_Z > 0$ and $A_{\alpha} < \infty$ by assumption,
this happens if 
\[ n^{1-\alpha}\sigma^{2\alpha}\ /\ [\log\set{1 + (\mu_1 n_{\sigma})^{-1}} + \log\set{N_{\poptreatment}(1/3n_{\sigma})}] \to \infty. \]
As the eigenvalue sum $\sum_{j}\mu_j$ is finite by assumption, each term $\mu_j$ is necessarily bounded by a constant multiple of the corresponding term $1/j$
in the (infinite-sum) harmonic series $\sum_j 1/j$; thus, without loss of generality, \smash{$\mu_j \le cj^{-1/p}$} for $p \le 1$.
As discussed in \citep[Equation 14]{fischer2017sobolev}, it follows that
\smash{$N_{\poptreatment}(1/3n_{\sigma}) \le c n_{\sigma}^{1/p}$}. Thus, the denominator is logarithmic in $n$ and $\sigma$,
and so long as the numerator increases polynomially, the ratio will tend to infinity. As the numerator
is $(n\sigma^{(2\alpha)/(1-\alpha)})^{1-\alpha}$ and $\alpha < 1$, this happens
when \smash{$\sigma^2 \ge n^{\beta}$} where $\beta \alpha/(1-\alpha) > -1$,
i.e., where \smash{$\beta > 1-\alpha^{-1}$.} 

We conclude by showing that the two terms in \eqref{eq:uniform-imbalance-bound} are $o(n^{-1/2})$.
The first term is if $n_{\sigma}^{1-\kappa_m}\sum_{j=1}^{\infty} \min(\mu_j, n_{\sigma}^{-1}) \to 0$.
As shown in the proof of Lemma~\ref{lemma:complexity} above, the sum involved goes to zero if $n_{\sigma}^{-1}$ does.
And because we have assumed $\kappa_m \ge 1$, it follows that the first term is negligible if $n_{\sigma} \to \infty$,
i.e., if $\sigma^2 \ll n$. 
The second term is negligible if $n_{\sigma}^{\kappa} \gg n$ for $\kappa = \kappa_m + \kappa_g + \kappa_g (1-\kappa_m)$.
As $n_\sigma^{\kappa}/n \sim n^{\kappa-1}/\sigma^{2\kappa}$, this happens if 
$n^{1-\kappa^{-1}} \gg \sigma^2$. 
Summarizing, bias is negligible if $\norm{m_0}_{\hh} < \infty$ and 
$\sigma$ satisfies \smash{$n^{1-\kappa^{-1}} \gg \sigma^2 \ge n^{1-\alpha^{-1}+\epsilon}$} for $\epsilon>0$.

\subsection{Convergence of weights}
\label{sec:convergence-of-weights}
Our bound on the difference between the estimated weights $\hriesz_i$ 
and the weights $\riesz[\psi](X_i,W_i)$ is based on a reduction to Theorem 2 of \citet{hirshberg2018augmented}.
In order to apply it, we first show the weights $\hgamma$ that we discuss here are an instance of the weights $\hgamma$ discussed in that paper.
\begin{proposition}
\label{prop:simplification-of-duality}
Let $h(x,w,f) = T(x,w)f(x,0)$, let $\bb$ be an absolutely convex set of functions on a set $\xx$,
and let $\bb^C$ be the unit ball for the cartesian product of $C+1$ copies of this space considered as functions 
$f(x,w)$ on $(\xx, \set{0 \ldots C})$. Then 
\begin{equation}
\label{eq:weights-minimax-linear-amle-form}
\ell_{n,\bb^C}(\riesz) = I_{h,\bb^C}^2(\riesz) + \lambda \norm{\riesz}^2,\ I_{h,\ff}=\sup_{f \in \ff}\frac{1}{n}\sum_{i=1}^n\{ h(X_i,W_i,f) - \riesz_i f(X_i,W_i) \},
\end{equation}
has a unique minimum at $\hriesz$ satisfying $\hriesz_i = 0$ for $W_i \neq 0$. 
\end{proposition}
Here we take $\bb$ to be the unit ball of our RKHS $\hh$. 
Subject to the constraint that $\hriesz_i=0$ if $W_i \neq 0$, a constraint that is satisfied by the solution of \eqref{eq:weights-minimax-linear-amle-form},
this problem reduces to the problem \eqref{eq:ridge-equivalence-weighting} that defines our weights. Thus, our weights solve it.
Having established this, \citet[Theorem 2]{hirshberg2018augmented} 
implies that the following bound holds with probability $1 - 4\delta - 3\exp(-cnr_Q^2)$. 
\begin{equation}
\label{eq:weight-consistency-abstract}
\begin{aligned}
&n^{-1}\sum_{i:W_i=1}\{\hgamma_i - \tilde\riesz(W_i,X_i)\}^2 \le 
6(nr^4/\sigma^2 + \norm{\tilde\riesz}_{\bb^C}r^2) \vee 8r^2  \quad \text{ for } \\
&\quad \tilde\riesz = \argmin_{\riesz} \norm{\riesz - \riesz[\psi]}_{L_2(\Pn)}^2 + (\sigma^2/n) \norm{\riesz}_{\bb^C}^2 \ \text{ and } \
r=r_Q \vee r_M \text{ where } \\
&\quad r_Q = \inf\set{ r > 0 : R_n(\bb_{cr}) \le cr^2 }, \\
&\quad r_M = \inf\set{ r > 0 : R_n\{(T-\riesz[\psi])\bb_r\} \le \delta r^2/2 }. \\
\end{aligned}
\end{equation}
Here $\norm{f}_{\bb^C}$ is the \emph{gauge} $\inf\set{ \alpha > 0 : f/\alpha \in \bb^C}$,
which is equal to $\max_w \norm{f(w,\cdot)}_{\hh}$.

We'll now simpify this. First, observe that $\tilde\riesz(w,x)=0$ for $w \neq 0$. 
It inherits this property from $\riesz[\psi](w,x)$, as changing $\riesz(w,x)$ for $w \neq 0$ increases
both $\norm{\riesz - \riesz[\psi]}_{L_2(\Pn)}$ and $\norm{\riesz}_{\bb^C}$. 
Thus, $\tilde \riesz(w,x) = \ind{w=0}\tilde g(x)$ where
\[ \tilde g = \argmin_g \ell_n(g),\ \ell_n(g) = n^{-1}\sum_{i=1}^n \ind{W_i=0} \{ g(X_i) - g_{\psi}(X_i)\}^2 + (\sigma^2/n) \norm{g}_{\hh}^2. \]
We will show that $\ell_n(\tilde g) = o_p(1)$ by considering the minimizer $g_{\delta}$ of 
\begin{align*}
\ell_\delta(g) &= \delta^{-1}\E[ \ind{W_i=0} \{g(X_i) - g_{\psi}(X_i)\}^2] + (\sigma^2/n) \norm{g}_{\hh}^2 \\
	       &= (p_Z/\delta)\norm{g-g_{\psi}}_{L_2(\poptreatment)}^2 + (\sigma^2/n) \norm{g}_{\hh}^2.
\end{align*}
By Markov's inequality, $\ell_n(g_\delta) \le \ell_{\delta}(g_\delta)$ on an event of probability $1-\delta$,
and because $\ell_n(\tilde g) \le \ell_n(g_\delta)$, on that event $\ell_n(\tilde g) \le \ell_{\delta}(g_\delta)$.
And because $\hh$ is dense in $L_2(\poptreatment)$ and $g_{\psi} \in L_2(\poptreatment)$, there exists
a sequence $g_n$ converging to $g_{\psi}$ with $\norm{g_n}_{\hh} \to \infty$ arbitrarily slowly.
In particular, whenever $\sigma^2 \ll n$,
there exists such a sequence with $\norm{g_n}_{\hh} \ll n/\sigma^2$,
so $\ell_{\delta}(g_{\delta}) \le \ell_{\delta}(g_n) \to 0$.
As $\norm{\tilde g}_{\hh} \le (n/\sigma^2)^{1/2}\ell_n(\tilde g)^{1/2}$, 
it follows that on an event of probability $1-\delta$, 
$\norm{\tilde g}_{\hh} \ll (n/\sigma^2)^{1/2}$,
and by the union bound that, on an event of probability $1- 5\delta - 3\exp(-cnr_Q^2)$ for $r$ and $r_Q$ as in \eqref{eq:weight-consistency-abstract},
\begin{equation*}
\label{eq:weight-consistency-abstract-2}
n^{-1}\sum_{i:W_i=0}\{\hgamma_i - \tilde\riesz(W_i,X_i)\}^2 
\le 6 \{ nr^4/\sigma^2 + o\{(n/\sigma^2)^{1/2}\} r^2\} \vee 8r^2.
\end{equation*}
This bound holds with high probability and converges to zero if $r \to 0$ and $n \gg \sigma^2 \gg nr^4$.
Substituting an upper bound on $r$ yields a sufficient condition, so we conclude by bounding $r=r_Q \vee r_M$.

For $r_{\star}$ as in Lemma~\ref{lemma:complexity}, comparing definitions, $r_Q=cr_{\star}$.
Furthermore, under our assumption that $g_{\psi}$ and therefore $\riesz[\psi]$ is bounded, 
$R_n\{(T-\riesz[\psi])\hh_r\} \le cR_n(\hh_{r})$ by the contraction inequality for Rademacher processes \citep[e.g.,][Theorem 3.1.21]{gine2015mathematical},
so $r_M \le cr_{\star}$. Thus, the rates for $r_{\star}$ from Lemma~\ref{lemma:complexity} imply that 
the lower bound in the sufficient condition $n \gg \sigma^2 \gg nr^4$ is $o(1)$, and if $\mu_j \le cj^{-1/p}$ for $p < 1$, 
it is $o(n^{1-1/p})$.

Under the same condition and on the same event, $n^{-1}\sum_{i:W_i=0}\{\hgamma_i - \riesz[\psi](W_i,X_i)\}^2 \to 0$.
As we have shown above that \smash{$\ell_n(\tilde g)$} and therefore 
\smash{$\norm{\tilde\riesz - \riesz[\psi]}_{L_2(\Pn)}$} converges to zero
on this event, this follows by the triangle inequality.
And as discussed in \citet{hirshberg2018augmented}, via Chebyshev's inequality conditional on 
$W_1,X_1 \ldots W_n,X_n$, this implies that the difference
$n^{-1}\sum_{i:W_i=0}\{\hgamma_i - \riesz[\psi](W_i,X_i)\}\varepsilon_i$ between 
the noise term in our error decomposition \eqref{eq:linear-error-decomp}
and the `oracle noise term' $n^{-1}\sum_{i:W_i=0}\riesz[\psi](W_i,X_i)\varepsilon_i$
is $o_p(n^{-1/2})$.

\subsection{Conclusion}

The results of the previous two subsections hold under Assumption~\ref{assu:sigma}.
Recalling the error decomposition \eqref{eq:linear-error-decomp}, these results imply that
\[ \hat \psi_{ML} - \psi(m) = \tilde \psi(m) - \psi(m) + n^{-1}\sum_{i:W_i=0}\riesz[\psi](W_i,X_i)\varepsilon_i + o_p(n^{-1/2}). \]
This is \eqref{eq:scaled-targeted-asymptotically-linear}, as the non-negligible terms on the right side above 
are $n^{-1}\sum_{i=1}^n \iota(W_i,X_i,Y_i)$.

We can derive \eqref{eq:targeted-asymptotically-linear} from this straightforwardly.
Observe that  
\[ n^{1/2}\{ \hat p_T^{-1}\hat{\psi}_{ML}  - \mu_0\} =
\hat p_T^{-1}n^{1/2}[\{\hat \psi_{ML} - P(T_i=1)\mu_0\} 
                  + \mu_0 \{ P(T_i=1) - \hat p_T \}].  
\]
If $n^{1/2}$ times the bracketed quantity converges to $Z_n$, the product converges to $Z_n/P(T_i=1)$.
To characterize $Z_n$, observe that via \eqref{eq:scaled-targeted-asymptotically-linear}, $n^{1/2}$ times the bracketed quantity converges to
\begin{align*}
&n^{-1/2}\sum_i \{ \influence(W_i,X_i,Y_i) + \mu_0 P(T_i=1) - T_i \mu_0 \}\ \\
&= n^{-1/2} \sum_i \influence''(W_i,X_i,Y_i),\ \influence''(w,x,y) = T(w,x) (m_0(x) - \mu_0) + \riesz[\psi](x,w)\{y - \mo(x)\}. 
\end{align*}
Thus, the influence function of $p_T^{-1}\hat{\psi}_{ML}$ is $\influence'(w,x,y) = \influence''(w,x,y)/P(T_i=1)$ 
as claimed in \eqref{eq:targeted-asymptotically-linear}.

\section{Proofs for lemmas used in Section~\ref{sec:proving-finite-sample-bounds}}
\label{sec:minimax-linear-lemmas}
Here we collect proofs for all the lemmas and propositions stated in the previous section.

\begin{proof}[\bkaproof of Lemma~\ref{lemma:ridge-equivalence}]

To simplify our notation, we'll use $\treatment_i$ as shorthand for $\ind{W_i=0}$. Our weights, defined in \eqref{eq:ridge-equivalence-weighting}, minimize
\begin{align*}
&\frac{\sigma^2}{n^2}\sum_{i : \treatment_i=1 } \gamma_i^2 + \sup_{f : \norm{f} \le 1}\{\frac{1}{n}\sum_{i} (\target_i - \treatment_i \gamma_i) \inner{K_{X_i}, f}\}^2 \\
&=\frac{\sigma^2}{n^2}\sum_{i : \treatment_i=1} \gamma_i^2 + \inner*{\frac{1}{n}\sum_i (\target_i - \treatment_i \gamma_i) K_{X_i}, 
																\frac{1}{n}\sum_j (\target_j - \treatment_j \gamma_j) K_{X_j}} \\
=&\frac{\sigma^2}{n^2}\sum_{i : \treatment_i=1} \gamma_i^2 + \frac{1}{n^2}\sum_{i,j} (\target_i - \treatment_i \gamma_i)(\target_j - \treatment_j \gamma_j) K(X_i,X_j) \\
=&\frac{1}{n^2}(\sigma^2 \gamma^T \gamma + 1^T K_{\target,\target} 1 - 2\gamma^T K_{\treatment,\target} 1 + \gamma^T K_{\treatment,\treatment} \gamma)  \\
=&\frac{1}{n^2}\{1^T K_{\target,\target} 1 - 2\gamma^T K_{\treatment,\target} 1 + \gamma^T \p{K_{\treatment,\treatment} + \sigma^2 I} \gamma\} 
\end{align*}
where $K$ is the Gram matrix ($K_{i,j} = K(X_i,X_j)$), $1$ is a vector of $\abs{\set{i: \target_i=1}}$ ones, and subscripting by $\treatment$ or $\target$ takes the rows of columns corresponding to units in those groups. At the minimum over $\gamma$, the derivative with respect to $\gamma$ will be zero, so our weights solve $(K_{\treatment,\treatment} + \sigma^2 I)\gamma = K_{\treatment,\target}1$, and the weighted average of treatment outcomes is
\begin{equation}
\label{eq:explicit-ridge-weighting}
 n^{-1}\sum_{i=1}^n\treatment_i \hgamma_i Y_i=n^{-1}Y_{\treatment}^T \hgamma = n^{-1}Y_{\treatment}^T (K_{\treatment,\treatment} + \sigma^2 I)^{-1} K_{\treatment,\target}1.
\end{equation}

Now consider ridge regression on the treated units. We estimate $\hmo$ solving 
\[ \min_{f} \sum_{i : \treatment_i=1} \p{Y_i - \inner{K_{X_i}, f}}^2 + \sigma^2 \norm{f}^2. \]
We can write it equivalently in constrained form,
\[ \min_{r,f} \sum_{i : \treatment_i=1} r_i^2 + \sigma^2 \norm{f}^2 \ \text{ where } \ r_i = \inner{K_{X_i}, f} - Y_i. \]
This problem is solved by a saddle point of the Lagrangian \citep[Theorem 3.6.8]{peypouquet2015convex},
\begin{align*}
&L((r,f), \lambda) = \sum_{i : \treatment_i=1} r_i^2 + \sigma^2 \norm{f}^2  + 2\sum_{i:\treatment_i=1} \lambda_i \p{\inner{K_{X_i}, f} - Y_i - r_i}. \\
\end{align*}
For given $\lambda$, we can minimize over $(r,f)$ explicitly, solving  the conditions
$r_i - \lambda_i = 0$ and $\sigma^2 f + \sum_{i : \treatment_i=1} \lambda_i K_{X_i} = 0$
that arise from setting the derivatives with respect to $r_i$ and $f$ to zero. 
Substituting the optimal values 
$\hat{r}_i = \lambda_i$ and $\hat{f} = -\sigma^{-2}\sum_{i : \treatment_i=1} \lambda_i K_{X_i}$,  
\begin{align*}
& L((\hat{r},\hat{f}),\lambda)) = \sum_{i \in \treatment} \lambda_i^2  + \sigma^{-2} \inner*{\sum_{i : \treatment_i=1}\lambda_i K{X_i},\sum_{j : \treatment_j=1} \lambda_j K_{X_j}}  \\
& \quad + 2\sum_{i : \treatment_i=1} \lambda_i \p{-\sigma^{-2}\inner*{\sum_{j : \treatment_j=1}\lambda_j K_{X_j} , K_{X_i}} - Y_i - \lambda_i}   \\
&=\sum_{i : \treatment_i=1} \lambda_i^2 + \sigma^{-2} \sum_{i,j : \treatment_i=\treatment_j=1}\lambda_i \lambda_j K(X_i,X_j)  \\ 
&	-2\sigma^{-2}\sum_{i,j : \treatment_i=\treatment_j=1}\lambda_i \lambda_j K(X_j,X_i) - 2\sum_{i : \treatment_i=1}\p{ \lambda_i Y_i + \lambda_i^2}   \\
&=-\lambda^T \lambda  - 2\lambda^T Y_{\treatment} - \sigma^{-2} \lambda^T K_{\treatment,\treatment} \lambda \\
&= -2\lambda^T Y_{\treatment} - \lambda^T\p{\sigma^{-2}K_{\treatment,\treatment} + I}\lambda.
\end{align*}
This is maximized at $\hat{\lambda} = -\p{\sigma^{-2}K_{\treatment,\treatment} + I}^{-1}Y_{\treatment} = -\sigma^2 \p{K_{\treatment,\treatment} + \sigma^2 I}^{-1} Y_{\treatment}$. 
Thus, we have a saddle at $((\hat{r}, \hat{f}), \hat{\lambda})$ and the function $\hmo$ solving our problem is $\hat{f}$.
Substituting in $\hat{\lambda}$ into our expression for $\hat{f}$ above,
\begin{align*}
\inner{K_{x}, \hat{f}} &= \inner*{-\sigma^{-2}\sum_{i : \treatment_i=1} \{-\sigma^2 \p{K_{\treatment,\treatment} + \sigma^2 I}^{-1}_{i,\treatment} Y_{\treatment}\} K_{X_i}, K_{x}} \\
&= \sum_{i : \treatment_i=1} Y_{\treatment}^T \p{K_{\treatment,\treatment} + \sigma^2 I}^{-1}_{\treatment,i} K(X_i,x). \\
\end{align*}
Therefore our ridge regression prediction $\hat{f}$, averaged over our target sample, is
\begin{align*}
\inner*{n^{-1}\sum_{j : \target_j=1} K_{X_j}, \hat{f}} 
&= n^{-1}\sum_{j : \target_j=1}\sum_{i : \treatment_i=1} Y_{\treatment}^T \p{K_{\treatment,\treatment} + \sigma^2 I}^{-1}_{\treatment,i} K(X_i,X_j). \\
&= n^{-1}Y_{\treatment}^T \p{K_{\treatment,\treatment} + \sigma^2 I}^{-1}K_{\treatment,\target} 1.
\end{align*}
This is the weighted average of treatment outcomes using our minimax weights, completing our proof.
\end{proof}

\begin{proof}[\bkaproof of Lemma~\ref{lemma:rkhs-approx-extremal}]
For $\eta=1/2$, the bound holds for $\tilde g = g$. For $\eta < 1/2$,
we reparameterize \citet[Theorem 4.1]{cucker2007learning} in terms of $\eta = \theta/(4+2\theta)$ to get the bound 
\[ \inf_{\norm{\tilde g} \le R} \norm{\tilde g - g} \le (2\norm{g}_{L_2(\eta)})^{\frac{1}{1-2\eta}} R^{-\frac{2\eta}{1-2\eta}}. \]
If this infimum is attained at some $\tilde g$, then 
\[ \norm{\tilde g - g}_{L_2(\eta)} + t \norm{\tilde g} \le (2\norm{g}_{L_2(\eta)})^{\frac{1}{1-2\eta}} R^{-\frac{2\eta}{1-2\eta}} + tR, \]
and the infimum over $\tilde g$ of $\norm{\tilde g - g}_{L_2(\eta)} + t \norm{\tilde g}$ can be no larger than the minimum of the right side above over $R$.
This bound is infinite at $R=0$ and strictly convex and differentiable for $R > 0$, so it is minimized 
over $R \ge 0$ at the zero of its derivative. This which occurs where
\[ t = (2\norm{g}_{L_2(\eta)})^{\frac{1}{1-2\eta}} \frac{2\eta}{1-2\eta} R^{-\p{1+\frac{2\eta}{1-2\eta}}} \ \text{ i.e. }
   R = t^{-(1-2\eta)} 2\norm{g}_{L_2(\eta)}\p{\frac{2\eta}{1-2\eta}}^{1-2\eta}. \]
Evaluating our bound at this value of $R$ yields our bound.
\end{proof}

\begin{proof}[\bkaproof of Lemma~\ref{lemma:regularized-product-opnorm-bound}]
The eigenvalues of the product have the form $\sigma^\nu / (\sigma + \lambda) = 1 / (\sigma^{1-\nu} + \lambda \sigma^{-\nu})$
where $\sigma$ is an eigenvalue of $L$. We maximize this expression over all $\sigma$. 
The derivative of the denominator is 
\[ (1-\nu)\sigma^{-\nu} - \lambda \nu \sigma^{-(1+\nu)} = \sigma^{-(1+\nu)} \{ (1-\nu)\sigma - \lambda \nu \} \]
and is zero at $\sigma = \lambda \nu / ( 1-\nu)$. To check that this is a minimum of the denominator, we check positivity of the second derivative at this point.
The second derivative of the denominator is
\[-\nu(1-\nu)\sigma^{-(1-\nu)} + \lambda \nu (1+\nu) \sigma^{-(2+\nu)} = \nu \sigma^{-(2+\nu)} \{ -(1-\nu)\sigma + \lambda (1+\nu) \} \]
and at  $\sigma = \lambda \nu / ( 1-\nu)$ this takes the value $\nu \sigma^{-(2+\nu)} \lambda > 0$. 
Evaluating the expression $1 / (\sigma^{1-\nu} + \lambda \sigma^{-\nu})$ at $\sigma = \lambda \nu / ( 1-\nu)$ yields our bound.
\end{proof}

\begin{proof}[\bkaproof of Proposition~\ref{prop:simplification-of-duality}]
Consider the decomposition of the expression maximized in $I_{h,\bb^C}$,
\[ \frac{1}{n}\sum_{i=1}^n\{h(X_i,W_i,f) - \riesz_i f(X_i,W_i)\} = \frac{1}{n}\sum_{i=1}^n\{T(X_i,W_i)-\ind{W_i=0}\riesz_i\}f(x,0) 
+ \frac{1}{n}\sum_{i:W_i \neq 0} \riesz_i f(x,w). \]
If there were nonzero weights $\riesz_i$ in the second sum,
functions $f(\cdot,1) \ldots f(\cdot,C)$ could be chosen from the (symmetric) unit ball $\bb$
that make the second term match the first in sign. It follows that the weights $\riesz'_i = \ind{W_i=0}\riesz_i$
satisfy \smash{$I_{h,\bb^C}^2(\riesz')$} \smash{$\le I_{h,\bb^C}^2(\riesz)$} and,  unless $\riesz_i'=\riesz_i$, 
\smash{$\norm{\riesz'}^2 < \norm{\riesz'}^2$}.
Therefore it suffices to optimize over weights of the form $\ind{W_i=0}\gamma_i$.
\end{proof}

\section{Additional Simulation Results}
\label{sec:additional-simulations}
In this section, we include simulation results for the classic example of \citet*{kang2007demystifying},
as well as variations on the example of \citet{hainmueller} discussed in Section~\ref{sec:empirical-performance}.

The Kang and Schafer example was designed to illustrate that methods using estimated inverse propensity weights can be unstable.
Here, the observations $X_i \in \R^4, W_i \in \set{0,1}, Y_i \in \R$ are defined in terms of a latent vector of standard normal random variables $Z_i \in \R^4$: $X_{i1} = \exp(Z_{i1}/2), X_{i2}=Z_{i2}/\{1+\exp(Z_{i1}) + 10\}, X_{i3} = (Z_{i1}Z_{i3}/25 + .06)^3, X_{i4} = (Z_{i2} + Z_{i4} + 20)^2$; $\pr(W_i=0 \mid Z_i) = \textrm{logit}^{-1}(-Z_{i1} + 0.5 Z_{i2} - 0.25Z_{i3} - 0.1 Z_{i4})$; and $Y_i = 210 + 27.4 Z_{i1} + 13.7 (Z_{12} + Z_{i3} + Z_{i4}) + \sigma_{\varepsilon} \varepsilon_i$ for a standard normal error $\varepsilon_i$ when $W_i=0$. 
In this example, the instability of the IPW and AIPW estimators persists even into large sample sizes, while the OLS estimator performs extremely well even in small samples. 
These phenomena are explained in detail by a comment on \citet*{kang2007demystifying} by \citet{robins2007comment}. 
In summary, there are regions of poor overlap between the distributions of the covariate $X_i$ between the treated and untreated subpopulations, which results in large inverse probability weights and therefore instability, but $\mo(x)$ is sufficiently linear throughout the support of $X_i$ that an estimator fit on the treated units extrapolates well into these regions of poor overlap. 
In Figure \ref{table:ks}, we show that our estimator $\hat{\psi}_{ML_t}$, while not reliant on the linearity of $\mo(x)$, also performs very well in all sample sizes.
Furthermore, when the sample size is small and the noise level $\sigma_{\varepsilon}$ is large, the regularization of our implicit estimator of $\hat{m}$ improves RMSE --- in these settings, $\hat{\psi}_{ML}$ and $\hat{\psi}_{ML_t}$ with larger values of the tuning parameter $\sigma$ outperform OLS.

In the example of Hainmueller discussed in Section~\ref{sec:empirical-performance}, 
we observe $X_i \in \R^6$ with $X_{i1} \ldots X_{i3}$ jointly normal with mean zero and covariance matrix $\Sigma$ defined below
$X_{i4} \sim \operatorname{Uniform}([-3,3])$, $X_{i5} \sim \chi^2_1$, and $X_{i6} \sim \operatorname{Bernoulli}(1/2)$ independent of each other and independent of $X_{i1} \ldots X_{i3}$; missingness follows a probit model $pr(W_i = 0 \mid X_i) = \Psi\{\eta^{-1}(X_{i1} + 2 X_{i2} - 2 X_{i3} - X_{i4} - 0.5 X_{i5} + X_{i6})\}$; and outcomes follow a quadratic model $Y_i = (X_{i1} + X_{i2} + X_{i5})^2 + \sigma_{\varepsilon} \varepsilon_i$. This is outcome design 3 from \citet{hainmueller}. Here we show additional results for that design,
varying our overlap and noise level parameters $\eta$ and $\sigma_{\varepsilon}$, 
as well as for outcome designs 1 and 2, in which we have $Y_i = X_{i1} + X_{i2} + X_{i3} - X_{i4} + X_{i5} + X_{i6} + \sigma_{\varepsilon}\varepsilon_i$
and $Y_i = X_{i1} + X_{i2} + 0.2 X_{i3} X_{i4} - \sqrt{X_{i5}} + \sigma_{\varepsilon}\varepsilon_i$ respectively. 
As these outcomes models are more linear than outcome design 3, results in these variations are more favorable to OLS.
The covariance matrix of $X_{i1} \ldots X_{i3}$ used in all variations is 
\[ \Sigma = \begin{pmatrix} 2 & 1 & -1 \\ 1 & 1 & -0.5 \\ -1 & -0.5 & 1 \end{pmatrix}. \]
In this example, there is better overlap than in that of Kang and Schafer and the logit missingness model used in the IPW and AIPW methods is barely misspecified, 
so the IPW and AIPW perform well in moderate and large samples. 

For comparison, in addition to the estimators discussed in Section~\ref{sec:empirical-performance},
we also include the non-translation-invariant minimax linear estimator $\psi_{ML}$.
For this and $\psi_{ML_t}$, we take $\model$ to be the unit ball of the RKHS associated with the Mat\'ern Kernel, $K(x,y) = k_{\nu}(\norm{x - y})$ 
for $k_{\nu}(x) = \frac{(\sqrt{2 \nu}x)^{\nu}}{2^{\nu - 1} \Gamma(\nu)} BK_{\nu}(\sqrt{2\nu}x))$ where $BK_{\nu}$ is a modified Bessel function of the second kind. 
The RKHS associated with this kernel is the Sobolev space $H^s$ for $s=d/2 + \nu$ \citep{schaback2011missing}. 
We take $\nu$ to be $3/2$ and the primary level of the parameter $\sigma$ in \eqref{eq:ridge-equivalence-weighting} to be $0.1$,
although we will display some additional results for $\sigma=1$ and $\sigma=10$. 
Calculation of the estimators is straightforward, amounting to the solution of a symmetric $n \times n$ linear system, 
as discussed in the Proof of Lemma~\ref{lemma:ridge-equivalence} in Appendix~\ref{sec:minimax-linear-lemmas}.

\begin{figure}[!h]
\begin{center}
{\tiny  \begin{tabular}{||c|c||c|c|c|c||c|c|c|c||} 
\hline 
 & n & 50 & 200 & 1000 & 4000 & 50 & 200 & 1000 & 4000  \\ 
\hline 
&&\multicolumn{4}{c||}{RMSE} &  \multicolumn{4}{c||}{Half-width}  \\ 
&&\multicolumn{4}{c||}{Bias} &  \multicolumn{4}{c||}{Coverage}  \\ 
\hline 
\multirow{14}{*}{\rotatebox[origin=c]{90}{\parbox[c]{3cm}{\centering $\sigma_{\varepsilon}=1$}}} 
&\multirow{2}{*}{\centering IPW}
 & 40.8 & 79.3 & 126.5 & 169.6 & 17.6 & 17.5 & 26.7 & 51.5 \\ 
 &  & -3.1 & 19 & 33.7 & 55 & 0.53 & 0.51 & 0.25 & 0.08 \\ 
\cline{2-10} 
&\multirow{2}{*}{\centering AIPW}
 & 8.1 & 14.9 & 51.4 & 92 & 15.8 & 17 & 26.6 & 51.5 \\ 
 &  & -1.4 & -5.4 & -11.9 & -23.7 & 0.96 & 0.98 & 0.95 & 0.66 \\ 
\cline{2-10} 
&\multirow{2}{*}{\centering OLS}
 & 6.8 & 3.3 & 1.7 & 1.2 & 14.8 & 8 & 3.7 & 1.9 \\ 
 &  & 0.1 & -0.5 & -0.7 & -0.9 & 0.97 & 0.98 & 0.98 & 0.91 \\ 
\cline{2-10} 
&\multirow{2}{*}{\centering ML}
 & 38.2 & 20 & 8.8 & 4.3 & 16.5 & 7.5 & 3.6 & 1.9 \\ 
 &  & -36.7 & -19.5 & -8.7 & -4.3 & 0.01 & 0 & 0 & 0 \\ 
\cline{2-10} 
&\multirow{2}{*}{\centering MLt}
 & 8.9 & 5.2 & 2.6 & 1.5 & 12.9 & 7.1 & 3.5 & 1.9 \\ 
 &  & -6.6 & -4.4 & -2.3 & -1.3 & 0.85 & 0.83 & 0.83 & 0.82 \\ 
\cline{2-10} 
&\multirow{2}{*}{\centering MLt $10\sigma$}
 & 9.8 & 6.4 & 3.6 & 2.2 & 12 & 6.4 & 3.1 & 1.7 \\ 
 &  & -7.7 & -5.7 & -3.4 & -2.1 & 0.76 & 0.6 & 0.44 & 0.28 \\ 
\cline{2-10} 
&\multirow{2}{*}{\centering MLt $100\sigma$}
 & 11.9 & 10.2 & 8.7 & 7.1 & 11 & 5.4 & 2.5 & 1.3 \\ 
 &  & -9.7 & -9.6 & -8.6 & -7.1 & 0.57 & 0.1 & 0 & 0 \\ 
\cline{2-10} 
\hline 
\multirow{14}{*}{\rotatebox[origin=c]{90}{\parbox[c]{3cm}{\centering $\sigma_{\varepsilon}=50$}}} 
&\multirow{2}{*}{\centering IPW}
 & 40.5 & 308 & 479.7 & 842.3 & 33.9 & 35.9 & 54 & 86.7 \\ 
 &  & -3.3 & 30.8 & 59.4 & 78.7 & 0.75 & 0.74 & 0.5 & 0.15 \\ 
\cline{2-10} 
&\multirow{2}{*}{\centering AIPW}
 & 15.4 & 69.8 & 143.4 & 513.4 & 32.6 & 35.2 & 53.8 & 86.7 \\ 
 &  & -1.2 & -8.5 & -20.2 & -35.5 & 0.95 & 0.98 & 0.98 & 0.9 \\ 
\cline{2-10} 
&\multirow{2}{*}{\centering OLS}
 & 14.1 & 6.8 & 3.1 & 1.7 & 34.4 & 17.6 & 8 & 4 \\ 
 &  & 0.2 & -0.8 & -0.9 & -0.9 & 0.97 & 0.99 & 0.99 & 0.98 \\ 
\cline{2-10} 
&\multirow{2}{*}{\centering ML}
 & 40.1 & 20.8 & 9.4 & 4.5 & 29.9 & 17.2 & 8.4 & 4.4 \\ 
 &  & -37.5 & -19.6 & -8.8 & -4.2 & 0.27 & 0.37 & 0.47 & 0.53 \\ 
\cline{2-10} 
&\multirow{2}{*}{\centering MLt}
 & 14.2 & 8 & 3.9 & 2.1 & 27.9 & 17 & 8.4 & 4.4 \\ 
 &  & -6.9 & -4.5 & -2.4 & -1.3 & 0.94 & 0.96 & 0.97 & 0.98 \\ 
\cline{2-10} 
&\multirow{2}{*}{\centering MLt $10\sigma$}
 & 14.3 & 8.4 & 4.4 & 2.5 & 22.9 & 14.1 & 7.2 & 3.8 \\ 
 &  & -8.1 & -5.8 & -3.5 & -2 & 0.88 & 0.91 & 0.9 & 0.89 \\ 
\cline{2-10} 
&\multirow{2}{*}{\centering MLt $100\sigma$}
 & 15.8 & 11.5 & 9.1 & 7.1 & 12.8 & 5.9 & 3.7 & 2.6 \\ 
 &  & -10.1 & -9.7 & -8.7 & -7 & 0.53 & 0.25 & 0.03 & 0 \\ 
\cline{2-10} 
\hline 
\multirow{14}{*}{\rotatebox[origin=c]{90}{\parbox[c]{3cm}{\centering $\sigma_{\varepsilon}=200$}}} 
&\multirow{2}{*}{\centering IPW}
 & 68.3 & 631.3 & 5869.4 & 442.3 & 115.8 & 87.3 & 159.2 & 92.7 \\ 
 &  & -1.3 & 39 & 236 & 58.7 & 0.94 & 0.94 & 0.86 & 0.66 \\ 
\cline{2-10} 
&\multirow{2}{*}{\centering AIPW}
 & 55.6 & 55.5 & 687 & 345.4 & 115.2 & 85.7 & 158 & 92.7 \\ 
 &  & 0.1 & -3.9 & -41.8 & -27.3 & 0.95 & 0.98 & 1 & 1 \\ 
\cline{2-10} 
&\multirow{2}{*}{\centering OLS}
 & 51.1 & 23.3 & 10.4 & 5.1 & 124.5 & 62.9 & 28.3 & 14.2 \\ 
 &  & 1.3 & -0.4 & -0.9 & -0.8 & 0.97 & 0.99 & 1 & 1 \\ 
\cline{2-10} 
&\multirow{2}{*}{\centering ML}
 & 53.6 & 29.9 & 14.3 & 6.9 & 101.1 & 61.1 & 30.8 & 16.1 \\ 
 &  & -36.5 & -19.4 & -8.8 & -4.2 & 0.93 & 0.96 & 0.98 & 0.99 \\ 
\cline{2-10} 
&\multirow{2}{*}{\centering MLt}
 & 44.7 & 24.2 & 11.7 & 5.7 & 100.5 & 61 & 30.8 & 16.1 \\ 
 &  & -6.6 & -4.2 & -2.5 & -1.3 & 0.96 & 0.98 & 0.99 & 0.99 \\ 
\cline{2-10} 
&\multirow{2}{*}{\centering MLt $10\sigma$}
 & 42.8 & 22.2 & 10.5 & 5.4 & 79.4 & 49.9 & 25.9 & 13.8 \\ 
 &  & -7.5 & -5.6 & -3.5 & -2.1 & 0.92 & 0.96 & 0.99 & 0.99 \\ 
\cline{2-10} 
&\multirow{2}{*}{\centering MLt $100\sigma$}
 & 42.3 & 22.5 & 12.7 & 8.4 & 29.1 & 11.4 & 11.2 & 9 \\ 
 &  & -9.9 & -9.5 & -8.8 & -7.1 & 0.49 & 0.4 & 0.59 & 0.66 \\ 
\cline{2-10} 
\hline 
\multirow{14}{*}{\rotatebox[origin=c]{90}{\parbox[c]{3cm}{\centering $\sigma_{\varepsilon}=1000$}}} 
&\multirow{2}{*}{\centering IPW}
 & 305.6 & 458.6 & 627.1 & 24863.4 & 585.9 & 436.6 & 352.7 & 912.8 \\ 
 &  & -0.6 & 30.4 & 50.9 & 928.5 & 0.96 & 0.97 & 0.97 & 0.97 \\ 
\cline{2-10} 
&\multirow{2}{*}{\centering AIPW}
 & 287 & 281.9 & 355 & 5577.7 & 584.8 & 436 & 352.5 & 912.2 \\ 
 &  & -1.5 & 2.3 & -19.9 & -142.7 & 0.94 & 0.99 & 1 & 1 \\ 
\cline{2-10} 
&\multirow{2}{*}{\centering OLS}
 & 257.4 & 115.8 & 50.7 & 26.2 & 619.7 & 312.9 & 140.3 & 70.3 \\ 
 &  & -3.3 & 3.4 & 1.5 & -0.4 & 0.98 & 0.99 & 1 & 0.99 \\ 
\cline{2-10} 
&\multirow{2}{*}{\centering ML}
 & 199.7 & 111.8 & 56.4 & 29.5 & 497.1 & 305.7 & 152.6 & 79.9 \\ 
 &  & -39.3 & -18.8 & -7.3 & -3.7 & 0.99 & 0.99 & 1 & 0.99 \\ 
\cline{2-10} 
&\multirow{2}{*}{\centering MLt}
 & 224.9 & 116.7 & 57.1 & 29.6 & 496.9 & 305.6 & 152.6 & 79.9 \\ 
 &  & -8.4 & -3.3 & -0.8 & -0.8 & 0.97 & 0.99 & 1 & 1 \\ 
\cline{2-10} 
&\multirow{2}{*}{\centering MLt $10\sigma$}
 & 211.5 & 105.1 & 50.1 & 26 & 392.3 & 248.9 & 128.5 & 68.8 \\ 
 &  & -8.5 & -2 & -1.8 & -1.4 & 0.93 & 0.98 & 0.99 & 0.99 \\ 
\cline{2-10} 
&\multirow{2}{*}{\centering MLt $100\sigma$}
 & 205.3 & 99.3 & 46.4 & 24.6 & 137.8 & 51.5 & 54.7 & 44.5 \\ 
 &  & -9.6 & -5.6 & -6.8 & -7.1 & 0.46 & 0.37 & 0.75 & 0.93 \\ 
\cline{2-10} 
\hline 
\end{tabular}
\\ 
}
\end{center}
\caption[Estimation error in the example of \citet*{kang2007demystifying}]{Root mean squared error (RMSE), bias, and confidence interval half-width and coverage over 1000 replications in the example of \citet*{kang2007demystifying}. Here we take the tuning parameter $\sigma$ to be $0.1$ in the estimators ML and MLt. The notation MLt $10\sigma$ and $100\sigma$ indicates the substitution of $1$ and $10$ respectively.}
\label{table:ks}
\end{figure}

\begin{figure}[!h]
\begin{center}
{\tiny  \begin{tabular}{||c|c|c||c|c|c|c||c|c|c|c||} 
\hline 
 &  & n & 50 & 200 & 1000 & 4000 & 50 & 200 & 1000 & 4000  \\ 
\hline 
&&&\multicolumn{4}{c||}{rmse} &  \multicolumn{4}{c||}{half-width} \\ 
&&&\multicolumn{4}{c||}{bias} &  \multicolumn{4}{c||}{coverage}  \\ 
\hline 
\multirow{12}{*}{\rotatebox[origin=c]{90}{\parbox[l]{3cm}{\centering $\eta = \sqrt{100}$ \  (high overlap) }}} 
&\multirow{6}{*}{\rotatebox[origin=c]{90}{\parbox[c]{3cm}{\centering $\sigma_{\varepsilon}=1.0$}}} 
&\multirow{2}{*}{\centering IPW}
 & 0.72 & 0.26 & 0.12 & 0.06 & 1 & 0.5 & 0.23 & 0.11 \\ 
 &  &  & -0.08 & -0.02 & 0 & 0 & 0.88 & 0.94 & 0.95 & 0.95 \\ 
\cline{3-11} 
&&\multirow{2}{*}{\centering AIPW}
 & 0.49 & 0.22 & 0.1 & 0.05 & 0.98 & 0.5 & 0.23 & 0.11 \\ 
 &  &  & -0.03 & -0.01 & 0 & 0 & 0.94 & 0.96 & 0.97 & 0.96 \\ 
\cline{3-11} 
&&\multirow{2}{*}{\centering OLS}
 & 0.48 & 0.22 & 0.1 & 0.05 & 0.99 & 0.5 & 0.22 & 0.11 \\ 
 &  &  & -0.02 & -0.01 & 0 & 0 & 0.96 & 0.97 & 0.97 & 0.96 \\ 
\cline{3-11} 
&&\multirow{2}{*}{\centering ML}
 & 0.56 & 0.31 & 0.14 & 0.07 & 0.9 & 0.48 & 0.23 & 0.12 \\ 
 &  &  & -0.38 & -0.22 & -0.1 & -0.05 & 0.87 & 0.88 & 0.88 & 0.88 \\ 
\cline{3-11} 
&&\multirow{2}{*}{\centering MLt}
 & 0.55 & 0.24 & 0.11 & 0.06 & 0.9 & 0.48 & 0.23 & 0.12 \\ 
 &  &  & 0.02 & -0.03 & -0.02 & -0.02 & 0.9 & 0.94 & 0.96 & 0.96 \\ 
\cline{3-11} 
&&\multirow{2}{*}{\centering MLt $10\sigma$}
 & 0.57 & 0.25 & 0.11 & 0.07 & 0.86 & 0.45 & 0.21 & 0.11 \\ 
 &  &  & 0.04 & -0.04 & -0.05 & -0.04 & 0.86 & 0.92 & 0.94 & 0.91 \\ 
\cline{3-11} 
&&\multirow{2}{*}{\centering MLt $100\sigma$}
 & 0.62 & 0.32 & 0.16 & 0.06 & 0.81 & 0.4 & 0.18 & 0.1 \\ 
 &  &  & 0.16 & 0.15 & 0.1 & -0.01 & 0.81 & 0.79 & 0.72 & 0.87 \\ 
\cline{3-11} 
\cline{2-11} 
&\multirow{6}{*}{\rotatebox[origin=c]{90}{\parbox[c]{3cm}{\centering $\sigma_{\varepsilon}=10.0$}}} 
&\multirow{2}{*}{\centering IPW}
 & 4.44 & 1.18 & 0.53 & 0.26 & 6.15 & 3.16 & 1.42 & 0.71 \\ 
 &  &  & -0.1 & -0.01 & -0.02 & 0 & 0.94 & 0.99 & 1 & 0.99 \\ 
\cline{3-11} 
&&\multirow{2}{*}{\centering AIPW}
 & 4.26 & 1.14 & 0.53 & 0.26 & 6.13 & 3.16 & 1.42 & 0.71 \\ 
 &  &  & -0.03 & 0 & -0.02 & 0 & 0.92 & 0.99 & 1 & 0.99 \\ 
\cline{3-11} 
&&\multirow{2}{*}{\centering OLS}
 & 2.67 & 1.04 & 0.5 & 0.25 & 5.88 & 3.01 & 1.37 & 0.68 \\ 
 &  &  & 0.03 & 0.02 & -0.02 & 0 & 0.95 & 1 & 1 & 0.99 \\ 
\cline{3-11} 
&&\multirow{2}{*}{\centering ML}
 & 1.76 & 0.94 & 0.55 & 0.28 & 4.23 & 2.68 & 1.42 & 0.76 \\ 
 &  &  & -0.35 & -0.2 & -0.12 & -0.06 & 0.98 & 0.99 & 1 & 0.99 \\ 
\cline{3-11} 
&&\multirow{2}{*}{\centering MLt}
 & 2.21 & 1.02 & 0.55 & 0.28 & 4.23 & 2.68 & 1.42 & 0.76 \\ 
 &  &  & 0.04 & 0 & -0.04 & -0.02 & 0.93 & 0.99 & 1 & 0.99 \\ 
\cline{3-11} 
&&\multirow{2}{*}{\centering MLt $10\sigma$}
 & 2.17 & 0.98 & 0.51 & 0.25 & 3.26 & 2.2 & 1.23 & 0.67 \\ 
 &  &  & 0.05 & -0.01 & -0.07 & -0.04 & 0.85 & 0.98 & 1 & 0.99 \\ 
\cline{3-11} 
&&\multirow{2}{*}{\centering MLt $100\sigma$}
 & 2.16 & 0.99 & 0.47 & 0.24 & 1.8 & 0.56 & 0.36 & 0.35 \\ 
 &  &  & 0.18 & 0.18 & 0.09 & -0.01 & 0.57 & 0.4 & 0.55 & 0.85 \\ 
\cline{3-11} 
\cline{2-11} 
&\multirow{6}{*}{\rotatebox[origin=c]{90}{\parbox[c]{3cm}{\centering $\sigma_{\varepsilon}=100.0$}}} 
&\multirow{2}{*}{\centering IPW}
 & 24.52 & 12.85 & 5.18 & 2.52 & 56.1 & 31.19 & 13.9 & 7.03 \\ 
 &  &  & -1.23 & 0.51 & 0.25 & -0.21 & 0.95 & 0.99 & 0.99 & 1 \\ 
\cline{3-11} 
&&\multirow{2}{*}{\centering AIPW}
 & 26.2 & 12.51 & 5.18 & 2.53 & 56.04 & 31.18 & 13.9 & 7.03 \\ 
 &  &  & -1.54 & 0.47 & 0.27 & -0.22 & 0.93 & 0.99 & 0.99 & 1 \\ 
\cline{3-11} 
&&\multirow{2}{*}{\centering OLS}
 & 25.11 & 11.81 & 5.02 & 2.45 & 57.98 & 29.93 & 13.44 & 6.78 \\ 
 &  &  & -1.47 & 0.6 & 0.28 & -0.12 & 0.96 & 0.99 & 0.99 & 1 \\ 
\cline{3-11} 
&&\multirow{2}{*}{\centering ML}
 & 15.69 & 10.28 & 5.19 & 2.67 & 41.02 & 26.69 & 13.92 & 7.49 \\ 
 &  &  & -1.51 & 0.28 & 0.06 & -0.18 & 0.97 & 0.98 & 0.99 & 0.99 \\ 
\cline{3-11} 
&&\multirow{2}{*}{\centering MLt}
 & 20.36 & 11.38 & 5.38 & 2.7 & 40.99 & 26.69 & 13.92 & 7.49 \\ 
 &  &  & -1.67 & 0.58 & 0.15 & -0.15 & 0.93 & 0.98 & 0.99 & 0.99 \\ 
\cline{3-11} 
&&\multirow{2}{*}{\centering MLt $10\sigma$}
 & 20.04 & 10.75 & 4.9 & 2.46 & 31.19 & 21.79 & 12.07 & 6.6 \\ 
 &  &  & -1.89 & 0.64 & 0.12 & -0.16 & 0.87 & 0.96 & 0.99 & 0.99 \\ 
\cline{3-11} 
&&\multirow{2}{*}{\centering MLt $100\sigma$}
 & 19.79 & 10.39 & 4.68 & 2.32 & 16.06 & 3.91 & 3.13 & 3.34 \\ 
 &  &  & -1.91 & 0.88 & 0.25 & -0.08 & 0.59 & 0.29 & 0.46 & 0.83 \\ 
\cline{3-11} 
\cline{2-11} 
\hline 
\multirow{12}{*}{\rotatebox[origin=c]{90}{\parbox[l]{3cm}{\centering $\eta = \sqrt{30}$ \  (low overlap) }}} 
&\multirow{6}{*}{\rotatebox[origin=c]{90}{\parbox[c]{3cm}{\centering $\sigma_{\varepsilon}=1.0$}}} 
&\multirow{2}{*}{\centering IPW}
 & 0.92 & 0.83 & 0.2 & 0.1 & 1.05 & 0.56 & 0.26 & 0.13 \\ 
 &  &  & -0.1 & -0.06 & 0 & 0 & 0.86 & 0.87 & 0.86 & 0.84 \\ 
\cline{3-11} 
&&\multirow{2}{*}{\centering AIPW}
 & 0.53 & 0.25 & 0.11 & 0.06 & 1.03 & 0.56 & 0.26 & 0.13 \\ 
 &  &  & 0.03 & 0.01 & 0 & 0 & 0.94 & 0.97 & 0.98 & 0.96 \\ 
\cline{3-11} 
&&\multirow{2}{*}{\centering OLS}
 & 0.49 & 0.24 & 0.1 & 0.06 & 1.05 & 0.52 & 0.23 & 0.12 \\ 
 &  &  & 0.02 & 0.01 & 0 & 0 & 0.96 & 0.96 & 0.98 & 0.96 \\ 
\cline{3-11} 
&&\multirow{2}{*}{\centering ML}
 & 0.53 & 0.32 & 0.16 & 0.09 & 0.89 & 0.48 & 0.24 & 0.13 \\ 
 &  &  & -0.35 & -0.23 & -0.13 & -0.07 & 0.91 & 0.86 & 0.85 & 0.84 \\ 
\cline{3-11} 
&&\multirow{2}{*}{\centering MLt}
 & 0.56 & 0.28 & 0.11 & 0.06 & 0.89 & 0.48 & 0.24 & 0.13 \\ 
 &  &  & 0.13 & 0.05 & 0.01 & 0 & 0.9 & 0.91 & 0.96 & 0.96 \\ 
\cline{3-11} 
&&\multirow{2}{*}{\centering MLt $10\sigma$}
 & 0.58 & 0.28 & 0.11 & 0.06 & 0.85 & 0.45 & 0.22 & 0.12 \\ 
 &  &  & 0.15 & 0.04 & -0.01 & -0.01 & 0.87 & 0.89 & 0.96 & 0.95 \\ 
\cline{3-11} 
&&\multirow{2}{*}{\centering MLt $100\sigma$}
 & 0.65 & 0.4 & 0.21 & 0.09 & 0.81 & 0.4 & 0.18 & 0.09 \\ 
 &  &  & 0.28 & 0.25 & 0.17 & 0.06 & 0.77 & 0.67 & 0.52 & 0.67 \\ 
\cline{3-11} 
\cline{2-11} 
&\multirow{6}{*}{\rotatebox[origin=c]{90}{\parbox[c]{3cm}{\centering $\sigma_{\varepsilon}=10.0$}}} 
&\multirow{2}{*}{\centering IPW}
 & 3.02 & 2.09 & 0.69 & 0.35 & 5.88 & 3.83 & 1.79 & 0.91 \\ 
 &  &  & -0.08 & -0.05 & 0.02 & 0 & 0.94 & 0.98 & 0.99 & 1 \\ 
\cline{3-11} 
&&\multirow{2}{*}{\centering AIPW}
 & 3.45 & 1.87 & 0.69 & 0.35 & 5.86 & 3.83 & 1.79 & 0.91 \\ 
 &  &  & 0.07 & -0.03 & 0.03 & 0 & 0.89 & 0.98 & 0.99 & 0.99 \\ 
\cline{3-11} 
&&\multirow{2}{*}{\centering OLS}
 & 3.14 & 1.3 & 0.53 & 0.28 & 6.76 & 3.41 & 1.54 & 0.77 \\ 
 &  &  & 0.03 & -0.05 & 0.03 & -0.01 & 0.94 & 0.99 & 0.99 & 1 \\ 
\cline{3-11} 
&&\multirow{2}{*}{\centering ML}
 & 1.7 & 1.06 & 0.57 & 0.36 & 4.16 & 2.74 & 1.58 & 0.94 \\ 
 &  &  & -0.4 & -0.29 & -0.1 & -0.07 & 0.98 & 1 & 1 & 0.99 \\ 
\cline{3-11} 
&&\multirow{2}{*}{\centering MLt}
 & 2.21 & 1.16 & 0.59 & 0.36 & 4.16 & 2.74 & 1.58 & 0.94 \\ 
 &  &  & 0.06 & -0.02 & 0.04 & 0 & 0.92 & 0.97 & 0.99 & 0.99 \\ 
\cline{3-11} 
&&\multirow{2}{*}{\centering MLt $10\sigma$}
 & 2.15 & 1.09 & 0.49 & 0.3 & 3.23 & 2.19 & 1.29 & 0.77 \\ 
 &  &  & 0.09 & -0.02 & 0.03 & -0.02 & 0.86 & 0.96 & 0.99 & 0.99 \\ 
\cline{3-11} 
&&\multirow{2}{*}{\centering MLt $100\sigma$}
 & 2.16 & 1.07 & 0.48 & 0.25 & 1.88 & 0.57 & 0.36 & 0.34 \\ 
 &  &  & 0.22 & 0.2 & 0.21 & 0.06 & 0.57 & 0.41 & 0.55 & 0.83 \\ 
\cline{3-11} 
\cline{2-11} 
&\multirow{6}{*}{\rotatebox[origin=c]{90}{\parbox[c]{3cm}{\centering $\sigma_{\varepsilon}=100.0$}}} 
&\multirow{2}{*}{\centering IPW}
 & 26.15 & 16.87 & 6.74 & 3.46 & 57.54 & 37.53 & 18.09 & 9.19 \\ 
 &  &  & -0.41 & 0.7 & -0.14 & -0.11 & 0.98 & 0.98 & 1 & 0.99 \\ 
\cline{3-11} 
&&\multirow{2}{*}{\centering AIPW}
 & 30.95 & 16.09 & 6.86 & 3.55 & 57.32 & 37.52 & 18.09 & 9.19 \\ 
 &  &  & 0.22 & 0.51 & -0.16 & -0.09 & 0.91 & 0.98 & 0.99 & 0.99 \\ 
\cline{3-11} 
&&\multirow{2}{*}{\centering OLS}
 & 28.91 & 12.48 & 5.58 & 2.66 & 66.44 & 34.36 & 15.28 & 7.64 \\ 
 &  &  & 0.06 & 0.15 & -0.03 & 0.01 & 0.98 & 0.98 & 0.99 & 1 \\ 
\cline{3-11} 
&&\multirow{2}{*}{\centering ML}
 & 14.83 & 9.77 & 5.7 & 3.25 & 40.58 & 27.47 & 15.74 & 9.27 \\ 
 &  &  & -0.68 & 0.07 & -0.38 & -0.23 & 0.99 & 0.99 & 0.99 & 1 \\ 
\cline{3-11} 
&&\multirow{2}{*}{\centering MLt}
 & 19.69 & 11.04 & 5.94 & 3.29 & 40.54 & 27.47 & 15.74 & 9.27 \\ 
 &  &  & -0.26 & 0.5 & -0.24 & -0.16 & 0.95 & 0.98 & 0.99 & 1 \\ 
\cline{3-11} 
&&\multirow{2}{*}{\centering MLt $10\sigma$}
 & 19.25 & 10.29 & 5.15 & 2.79 & 31.09 & 21.78 & 12.82 & 7.57 \\ 
 &  &  & -0.15 & 0.71 & -0.03 & -0.13 & 0.9 & 0.96 & 0.99 & 1 \\ 
\cline{3-11} 
&&\multirow{2}{*}{\centering MLt $100\sigma$}
 & 19.13 & 9.93 & 4.5 & 2.26 & 16.37 & 4.21 & 3.11 & 3.27 \\ 
 &  &  & -0.05 & 1.3 & 0.3 & 0.05 & 0.57 & 0.32 & 0.53 & 0.85 \\ 
\cline{3-11} 
\cline{2-11} 
\hline 
\end{tabular}
\\ 
}
\end{center}
\caption[Estimation error in the example of \citet*{hainmueller} with Outcome Design 1]{Root mean squared error (rmse), bias, and confidence interval half-width and coverage over 1000 replications
of Outcome Design 1 from \citet*{hainmueller}.
Here we take the tuning parameter $\sigma$ to be $0.1$ in the estimators ML and MLt. The notation MLt $10\sigma$ and $100\sigma$ indicates the substitution of $1$ and $10$ respectively.}
\end{figure}

\begin{figure}[!h]
\begin{center}
{\tiny  \begin{tabular}{||c|c|c||c|c|c|c||c|c|c|c||} 
\hline 
 &  & n & 50 & 200 & 1000 & 4000 & 50 & 200 & 1000 & 4000  \\ 
\hline 
&&&\multicolumn{4}{c||}{rmse} &  \multicolumn{4}{c||}{half-width} \\ 
&&&\multicolumn{4}{c||}{bias} &  \multicolumn{4}{c||}{coverage}  \\ 
\hline 
\multirow{12}{*}{\rotatebox[origin=c]{90}{\parbox[l]{3cm}{\centering $\eta = \sqrt{100}$ \  (high overlap) }}} 
&\multirow{6}{*}{\rotatebox[origin=c]{90}{\parbox[c]{3cm}{\centering $\sigma_{\varepsilon}=1.0$}}} 
&\multirow{2}{*}{\centering IPW}
 & 0.66 & 0.28 & 0.12 & 0.06 & 0.9 & 0.48 & 0.21 & 0.11 \\ 
 &  &  & 0.13 & 0.03 & 0.02 & 0.01 & 0.87 & 0.92 & 0.92 & 0.94 \\ 
\cline{3-11} 
&&\multirow{2}{*}{\centering AIPW}
 & 0.45 & 0.2 & 0.09 & 0.05 & 0.89 & 0.48 & 0.21 & 0.11 \\ 
 &  &  & 0 & -0.01 & 0 & 0 & 0.94 & 0.98 & 0.98 & 0.99 \\ 
\cline{3-11} 
&&\multirow{2}{*}{\centering OLS}
 & 0.44 & 0.2 & 0.09 & 0.05 & 0.9 & 0.46 & 0.21 & 0.1 \\ 
 &  &  & -0.01 & -0.04 & -0.03 & -0.03 & 0.94 & 0.98 & 0.97 & 0.95 \\ 
\cline{3-11} 
&&\multirow{2}{*}{\centering ML}
 & 0.7 & 0.41 & 0.21 & 0.1 & 0.8 & 0.43 & 0.21 & 0.11 \\ 
 &  &  & 0.6 & 0.37 & 0.19 & 0.09 & 0.68 & 0.65 & 0.6 & 0.67 \\ 
\cline{3-11} 
&&\multirow{2}{*}{\centering MLt}
 & 0.77 & 0.43 & 0.21 & 0.1 & 0.8 & 0.43 & 0.21 & 0.11 \\ 
 &  &  & 0.62 & 0.37 & 0.19 & 0.09 & 0.65 & 0.61 & 0.6 & 0.66 \\ 
\cline{3-11} 
&&\multirow{2}{*}{\centering MLt $10\sigma$}
 & 0.85 & 0.54 & 0.3 & 0.15 & 0.75 & 0.4 & 0.2 & 0.1 \\ 
 &  &  & 0.71 & 0.5 & 0.28 & 0.15 & 0.52 & 0.34 & 0.17 & 0.17 \\ 
\cline{3-11} 
&&\multirow{2}{*}{\centering MLt $100\sigma$}
 & 0.96 & 0.84 & 0.78 & 0.65 & 0.71 & 0.35 & 0.16 & 0.08 \\ 
 &  &  & 0.83 & 0.8 & 0.77 & 0.65 & 0.41 & 0.03 & 0 & 0 \\ 
\cline{3-11} 
\cline{2-11} 
&\multirow{6}{*}{\rotatebox[origin=c]{90}{\parbox[c]{3cm}{\centering $\sigma_{\varepsilon}=10.0$}}} 
&\multirow{2}{*}{\centering IPW}
 & 3.17 & 1.2 & 0.52 & 0.25 & 5.9 & 3.09 & 1.42 & 0.71 \\ 
 &  &  & 0.03 & 0.01 & 0.02 & 0.02 & 0.94 & 0.99 & 0.99 & 1 \\ 
\cline{3-11} 
&&\multirow{2}{*}{\centering AIPW}
 & 2.86 & 1.17 & 0.52 & 0.24 & 5.89 & 3.09 & 1.42 & 0.71 \\ 
 &  &  & -0.03 & -0.02 & 0.01 & 0.01 & 0.95 & 0.99 & 0.99 & 0.99 \\ 
\cline{3-11} 
&&\multirow{2}{*}{\centering OLS}
 & 2.54 & 1.12 & 0.5 & 0.24 & 5.96 & 3.02 & 1.37 & 0.68 \\ 
 &  &  & -0.07 & -0.04 & -0.02 & -0.02 & 0.97 & 0.99 & 0.99 & 0.99 \\ 
\cline{3-11} 
&&\multirow{2}{*}{\centering ML}
 & 1.74 & 1.07 & 0.56 & 0.27 & 4.2 & 2.69 & 1.41 & 0.75 \\ 
 &  &  & 0.61 & 0.36 & 0.19 & 0.1 & 0.96 & 0.99 & 0.99 & 0.99 \\ 
\cline{3-11} 
&&\multirow{2}{*}{\centering MLt}
 & 2.19 & 1.17 & 0.58 & 0.28 & 4.2 & 2.69 & 1.41 & 0.75 \\ 
 &  &  & 0.61 & 0.36 & 0.2 & 0.1 & 0.92 & 0.98 & 0.99 & 0.99 \\ 
\cline{3-11} 
&&\multirow{2}{*}{\centering MLt $10\sigma$}
 & 2.18 & 1.16 & 0.57 & 0.28 & 3.21 & 2.21 & 1.22 & 0.66 \\ 
 &  &  & 0.71 & 0.48 & 0.29 & 0.15 & 0.84 & 0.93 & 0.97 & 0.99 \\ 
\cline{3-11} 
&&\multirow{2}{*}{\centering MLt $100\sigma$}
 & 2.21 & 1.29 & 0.91 & 0.69 & 1.73 & 0.52 & 0.35 & 0.34 \\ 
 &  &  & 0.82 & 0.79 & 0.78 & 0.66 & 0.54 & 0.3 & 0.16 & 0.08 \\ 
\cline{3-11} 
\cline{2-11} 
&\multirow{6}{*}{\rotatebox[origin=c]{90}{\parbox[c]{3cm}{\centering $\sigma_{\varepsilon}=100.0$}}} 
&\multirow{2}{*}{\centering IPW}
 & 27.55 & 12.42 & 5.08 & 2.62 & 56.77 & 30.81 & 14.11 & 7.01 \\ 
 &  &  & 1.38 & 0.61 & 0.12 & 0.01 & 0.96 & 0.99 & 0.99 & 0.99 \\ 
\cline{3-11} 
&&\multirow{2}{*}{\centering AIPW}
 & 27.88 & 12.42 & 5.1 & 2.63 & 56.69 & 30.81 & 14.11 & 7.01 \\ 
 &  &  & 0.99 & 0.56 & 0.12 & -0.01 & 0.94 & 0.98 & 0.99 & 0.99 \\ 
\cline{3-11} 
&&\multirow{2}{*}{\centering OLS}
 & 25.73 & 11.78 & 4.88 & 2.55 & 58.52 & 29.94 & 13.59 & 6.77 \\ 
 &  &  & 0.69 & 0.47 & 0.07 & -0.03 & 0.97 & 0.99 & 0.99 & 0.99 \\ 
\cline{3-11} 
&&\multirow{2}{*}{\centering ML}
 & 15.91 & 10.2 & 5.05 & 2.78 & 41.57 & 26.69 & 14.06 & 7.49 \\ 
 &  &  & 0.8 & 0.63 & 0.23 & 0.05 & 0.98 & 0.99 & 0.99 & 0.99 \\ 
\cline{3-11} 
&&\multirow{2}{*}{\centering MLt}
 & 20.34 & 11.32 & 5.23 & 2.82 & 41.55 & 26.69 & 14.06 & 7.49 \\ 
 &  &  & 0.95 & 0.65 & 0.24 & 0.06 & 0.95 & 0.98 & 0.99 & 0.99 \\ 
\cline{3-11} 
&&\multirow{2}{*}{\centering MLt $10\sigma$}
 & 19.86 & 10.85 & 4.76 & 2.57 & 31.63 & 21.82 & 12.17 & 6.6 \\ 
 &  &  & 1.28 & 0.71 & 0.31 & 0.15 & 0.88 & 0.96 & 0.99 & 0.99 \\ 
\cline{3-11} 
&&\multirow{2}{*}{\centering MLt $100\sigma$}
 & 19.58 & 10.42 & 4.49 & 2.45 & 16.32 & 4 & 3.14 & 3.34 \\ 
 &  &  & 1.48 & 1 & 0.83 & 0.67 & 0.57 & 0.29 & 0.52 & 0.83 \\ 
\cline{3-11} 
\cline{2-11} 
\hline 
\multirow{12}{*}{\rotatebox[origin=c]{90}{\parbox[l]{3cm}{\centering $\eta = \sqrt{30}$ \  (low overlap) }}} 
&\multirow{6}{*}{\rotatebox[origin=c]{90}{\parbox[c]{3cm}{\centering $\sigma_{\varepsilon}=1.0$}}} 
&\multirow{2}{*}{\centering IPW}
 & 1.37 & 0.8 & 0.24 & 0.15 & 1 & 0.58 & 0.26 & 0.13 \\ 
 &  &  & 0.21 & 0.09 & 0.09 & 0.1 & 0.75 & 0.74 & 0.7 & 0.53 \\ 
\cline{3-11} 
&&\multirow{2}{*}{\centering AIPW}
 & 0.52 & 0.29 & 0.11 & 0.06 & 0.96 & 0.57 & 0.26 & 0.13 \\ 
 &  &  & -0.05 & -0.02 & -0.01 & -0.01 & 0.92 & 0.96 & 0.97 & 0.97 \\ 
\cline{3-11} 
&&\multirow{2}{*}{\centering OLS}
 & 0.47 & 0.23 & 0.12 & 0.08 & 0.97 & 0.49 & 0.22 & 0.11 \\ 
 &  &  & -0.07 & -0.07 & -0.07 & -0.07 & 0.95 & 0.97 & 0.94 & 0.82 \\ 
\cline{3-11} 
&&\multirow{2}{*}{\centering ML}
 & 0.86 & 0.59 & 0.34 & 0.2 & 0.82 & 0.45 & 0.23 & 0.13 \\ 
 &  &  & 0.8 & 0.56 & 0.33 & 0.19 & 0.52 & 0.27 & 0.13 & 0.11 \\ 
\cline{3-11} 
&&\multirow{2}{*}{\centering MLt}
 & 1.01 & 0.66 & 0.37 & 0.21 & 0.83 & 0.45 & 0.23 & 0.13 \\ 
 &  &  & 0.91 & 0.62 & 0.36 & 0.21 & 0.44 & 0.21 & 0.1 & 0.07 \\ 
\cline{3-11} 
&&\multirow{2}{*}{\centering MLt $10\sigma$}
 & 1.12 & 0.81 & 0.49 & 0.3 & 0.79 & 0.42 & 0.2 & 0.11 \\ 
 &  &  & 1.03 & 0.78 & 0.48 & 0.29 & 0.29 & 0.03 & 0 & 0 \\ 
\cline{3-11} 
&&\multirow{2}{*}{\centering MLt $100\sigma$}
 & 1.27 & 1.19 & 1.11 & 0.95 & 0.76 & 0.37 & 0.17 & 0.09 \\ 
 &  &  & 1.18 & 1.17 & 1.11 & 0.95 & 0.15 & 0 & 0 & 0 \\ 
\cline{3-11} 
\cline{2-11} 
&\multirow{6}{*}{\rotatebox[origin=c]{90}{\parbox[c]{3cm}{\centering $\sigma_{\varepsilon}=10.0$}}} 
&\multirow{2}{*}{\centering IPW}
 & 4.37 & 2.2 & 0.7 & 0.37 & 6.23 & 3.83 & 1.76 & 0.92 \\ 
 &  &  & 0.14 & 0.03 & 0.11 & 0.09 & 0.95 & 0.98 & 0.98 & 0.99 \\ 
\cline{3-11} 
&&\multirow{2}{*}{\centering AIPW}
 & 3.61 & 1.88 & 0.68 & 0.36 & 6.19 & 3.83 & 1.76 & 0.92 \\ 
 &  &  & -0.13 & -0.07 & -0.01 & -0.02 & 0.89 & 0.97 & 0.99 & 0.99 \\ 
\cline{3-11} 
&&\multirow{2}{*}{\centering OLS}
 & 2.98 & 1.3 & 0.56 & 0.29 & 6.65 & 3.4 & 1.53 & 0.77 \\ 
 &  &  & -0.14 & -0.09 & -0.05 & -0.08 & 0.97 & 0.99 & 0.99 & 0.99 \\ 
\cline{3-11} 
&&\multirow{2}{*}{\centering ML}
 & 1.78 & 1.18 & 0.67 & 0.38 & 4.13 & 2.73 & 1.57 & 0.94 \\ 
 &  &  & 0.76 & 0.56 & 0.34 & 0.18 & 0.97 & 0.98 & 0.98 & 0.98 \\ 
\cline{3-11} 
&&\multirow{2}{*}{\centering MLt}
 & 2.3 & 1.33 & 0.7 & 0.39 & 4.13 & 2.73 & 1.57 & 0.94 \\ 
 &  &  & 0.84 & 0.61 & 0.37 & 0.19 & 0.9 & 0.94 & 0.97 & 0.98 \\ 
\cline{3-11} 
&&\multirow{2}{*}{\centering MLt $10\sigma$}
 & 2.3 & 1.32 & 0.72 & 0.41 & 3.21 & 2.18 & 1.29 & 0.76 \\ 
 &  &  & 0.95 & 0.75 & 0.5 & 0.28 & 0.82 & 0.88 & 0.92 & 0.95 \\ 
\cline{3-11} 
&&\multirow{2}{*}{\centering MLt $100\sigma$}
 & 2.33 & 1.53 & 1.23 & 0.97 & 1.85 & 0.54 & 0.35 & 0.34 \\ 
 &  &  & 1.09 & 1.12 & 1.14 & 0.94 & 0.54 & 0.24 & 0.04 & 0 \\ 
\cline{3-11} 
\cline{2-11} 
&\multirow{6}{*}{\rotatebox[origin=c]{90}{\parbox[c]{3cm}{\centering $\sigma_{\varepsilon}=100.0$}}} 
&\multirow{2}{*}{\centering IPW}
 & 34.3 & 15.81 & 7 & 3.33 & 60.44 & 37.87 & 17.89 & 9.23 \\ 
 &  &  & -1.09 & -0.24 & 0.15 & 0.31 & 0.96 & 0.98 & 0.99 & 0.99 \\ 
\cline{3-11} 
&&\multirow{2}{*}{\centering AIPW}
 & 35.79 & 16.24 & 7.19 & 3.41 & 60.22 & 37.87 & 17.89 & 9.23 \\ 
 &  &  & -1.76 & -0.39 & 0.06 & 0.2 & 0.9 & 0.98 & 0.99 & 0.99 \\ 
\cline{3-11} 
&&\multirow{2}{*}{\centering OLS}
 & 29.06 & 12.55 & 5.52 & 2.74 & 65.9 & 33.7 & 15.25 & 7.64 \\ 
 &  &  & -1.51 & -0.27 & 0.07 & 0.08 & 0.97 & 0.99 & 1 & 1 \\ 
\cline{3-11} 
&&\multirow{2}{*}{\centering ML}
 & 15.19 & 9.86 & 5.9 & 3.29 & 40.53 & 27.07 & 15.69 & 9.3 \\ 
 &  &  & 0.19 & 0.47 & 0.49 & 0.36 & 0.98 & 0.99 & 0.99 & 1 \\ 
\cline{3-11} 
&&\multirow{2}{*}{\centering MLt}
 & 20.08 & 11.14 & 6.16 & 3.34 & 40.5 & 27.07 & 15.69 & 9.3 \\ 
 &  &  & 0.08 & 0.53 & 0.52 & 0.38 & 0.94 & 0.98 & 0.99 & 1 \\ 
\cline{3-11} 
&&\multirow{2}{*}{\centering MLt $10\sigma$}
 & 19.47 & 10.39 & 5.33 & 2.84 & 31.08 & 21.47 & 12.82 & 7.57 \\ 
 &  &  & 0.3 & 0.69 & 0.63 & 0.49 & 0.87 & 0.95 & 0.98 & 0.99 \\ 
\cline{3-11} 
&&\multirow{2}{*}{\centering MLt $100\sigma$}
 & 19.16 & 9.84 & 4.69 & 2.5 & 16.85 & 4.14 & 3.11 & 3.27 \\ 
 &  &  & 0.6 & 1.14 & 1.2 & 1.09 & 0.61 & 0.31 & 0.5 & 0.81 \\ 
\cline{3-11} 
\cline{2-11} 
\hline 
\end{tabular}
\\ 
}
\end{center}
\caption[Estimation error in the example of \citet*{hainmueller} with Outcome Design 2]{Root mean squared error (rmse), bias, and confidence interval half-width and coverage over 1000 replications
of Outcome Design 2 from \citet*{hainmueller}.
Here we take the tuning parameter $\sigma$ to be $0.1$ in the estimators ML and MLt. The notation MLt $10\sigma$ and $100\sigma$ indicates the substitution of $1$ and $10$ respectively.}
\end{figure}

\begin{figure}[!h]
\begin{center}
{\tiny  \begin{tabular}{||c|c|c||c|c|c|c||c|c|c|c||} 
\hline 
 &  & n & 50 & 200 & 1000 & 4000 & 50 & 200 & 1000 & 4000  \\ 
\hline 
&&&\multicolumn{4}{c||}{RMSE} &  \multicolumn{4}{c||}{Half-width} \\ 
&&&\multicolumn{4}{c||}{Bias} &  \multicolumn{4}{c||}{Coverage}  \\ 
\hline 
\multirow{12}{*}{\rotatebox[origin=c]{90}{\parbox[l]{3cm}{\centering $\eta = \sqrt{100}$ \  (high overlap) }}} 
&\multirow{6}{*}{\rotatebox[origin=c]{90}{\parbox[c]{3cm}{\centering $\sigma_{\varepsilon}=1.0$}}} 
&\multirow{2}{*}{\centering IPW}
 & 4.17 & 1.84 & 0.77 & 0.34 & 6.9 & 4.41 & 2.25 & 1.15 \\ 
 &  &  & -0.49 & -0.18 & -0.08 & -0.06 & 0.89 & 0.97 & 0.99 & 1 \\ 
\cline{3-11} 
&&\multirow{2}{*}{\centering AIPW}
 & 4.73 & 1.82 & 0.85 & 0.4 & 6.88 & 4.41 & 2.25 & 1.15 \\ 
 &  &  & -0.54 & -0.28 & -0.11 & -0.09 & 0.84 & 0.95 & 0.99 & 0.99 \\ 
\cline{3-11} 
&&\multirow{2}{*}{\centering OLS}
 & 3.19 & 1.84 & 1.34 & 1.17 & 5.69 & 3.45 & 1.68 & 0.86 \\ 
 &  &  & -1.48 & -1.26 & -1.17 & -1.13 & 0.84 & 0.88 & 0.72 & 0.23 \\ 
\cline{3-11} 
&&\multirow{2}{*}{\centering ML}
 & 3.03 & 1.95 & 1.05 & 0.57 & 4.41 & 2.91 & 1.71 & 1.02 \\ 
 &  &  & -2.44 & -1.71 & -0.96 & -0.53 & 0.7 & 0.79 & 0.9 & 0.97 \\ 
\cline{3-11} 
&&\multirow{2}{*}{\centering MLt}
 & 2.43 & 1.24 & 0.7 & 0.41 & 4.41 & 2.91 & 1.71 & 1.02 \\ 
 &  &  & -0.31 & -0.64 & -0.53 & -0.35 & 0.87 & 0.94 & 0.99 & 1 \\ 
\cline{3-11} 
&&\multirow{2}{*}{\centering MLt $10\sigma$}
 & 2.51 & 1.36 & 1.01 & 0.68 & 3.86 & 2.42 & 1.42 & 0.86 \\ 
 &  &  & -0.2 & -0.89 & -0.92 & -0.65 & 0.83 & 0.87 & 0.83 & 0.81 \\ 
\cline{3-11} 
&&\multirow{2}{*}{\centering MLt $100\sigma$}
 & 3.22 & 1.73 & 0.64 & 0.8 & 3.44 & 1.76 & 0.82 & 0.46 \\ 
 &  &  & 0.94 & 0.86 & 0.24 & -0.77 & 0.75 & 0.71 & 0.83 & 0.1 \\ 
\cline{3-11} 
\cline{2-11} 
&\multirow{6}{*}{\rotatebox[origin=c]{90}{\parbox[c]{3cm}{\centering $\sigma_{\varepsilon}=10.0$}}} 
&\multirow{2}{*}{\centering IPW}
 & 8.03 & 1.91 & 0.89 & 0.43 & 9.96 & 5.43 & 2.63 & 1.35 \\ 
 &  &  & -0.16 & -0.18 & -0.04 & -0.05 & 0.96 & 0.99 & 1 & 1 \\ 
\cline{3-11} 
&&\multirow{2}{*}{\centering AIPW}
 & 12.66 & 2.18 & 0.99 & 0.48 & 9.95 & 5.43 & 2.63 & 1.35 \\ 
 &  &  & 0.03 & -0.17 & -0.08 & -0.07 & 0.94 & 0.97 & 0.99 & 0.99 \\ 
\cline{3-11} 
&&\multirow{2}{*}{\centering OLS}
 & 3.93 & 2.13 & 1.38 & 1.18 & 8.42 & 4.59 & 2.17 & 1.09 \\ 
 &  &  & -1.17 & -1.18 & -1.09 & -1.1 & 0.96 & 0.95 & 0.87 & 0.49 \\ 
\cline{3-11} 
&&\multirow{2}{*}{\centering ML}
 & 3.37 & 2.18 & 1.16 & 0.62 & 6.21 & 3.94 & 2.21 & 1.27 \\ 
 &  &  & -2.34 & -1.71 & -0.92 & -0.51 & 0.89 & 0.91 & 0.96 & 0.98 \\ 
\cline{3-11} 
&&\multirow{2}{*}{\centering MLt}
 & 3.24 & 1.65 & 0.87 & 0.49 & 6.21 & 3.94 & 2.21 & 1.27 \\ 
 &  &  & -0.17 & -0.66 & -0.49 & -0.33 & 0.93 & 0.97 & 0.99 & 0.99 \\ 
\cline{3-11} 
&&\multirow{2}{*}{\centering MLt $10\sigma$}
 & 3.27 & 1.73 & 1.11 & 0.71 & 5.13 & 3.25 & 1.87 & 1.09 \\ 
 &  &  & -0.07 & -0.91 & -0.89 & -0.63 & 0.87 & 0.91 & 0.91 & 0.91 \\ 
\cline{3-11} 
&&\multirow{2}{*}{\centering MLt $100\sigma$}
 & 3.93 & 1.97 & 0.81 & 0.81 & 3.91 & 1.79 & 0.89 & 0.57 \\ 
 &  &  & 1.08 & 0.81 & 0.28 & -0.74 & 0.7 & 0.61 & 0.73 & 0.31 \\ 
\cline{3-11} 
\cline{2-11} 
&\multirow{6}{*}{\rotatebox[origin=c]{90}{\parbox[c]{3cm}{\centering $\sigma_{\varepsilon}=100.0$}}} 
&\multirow{2}{*}{\centering IPW}
 & 27.56 & 12.24 & 5.07 & 2.61 & 56.72 & 31.34 & 14.22 & 7.13 \\ 
 &  &  & 1.56 & 0.3 & -0.16 & 0.03 & 0.96 & 0.99 & 0.99 & 0.99 \\ 
\cline{3-11} 
&&\multirow{2}{*}{\centering AIPW}
 & 26.88 & 12.37 & 5.12 & 2.61 & 56.64 & 31.34 & 14.22 & 7.13 \\ 
 &  &  & 1.78 & 0.3 & -0.2 & 0.01 & 0.95 & 0.99 & 0.99 & 0.99 \\ 
\cline{3-11} 
&&\multirow{2}{*}{\centering OLS}
 & 24.47 & 11.57 & 5.06 & 2.73 & 58.67 & 30.26 & 13.65 & 6.84 \\ 
 &  &  & 0.71 & -0.9 & -1.28 & -1.04 & 0.98 & 1 & 0.99 & 0.99 \\ 
\cline{3-11} 
&&\multirow{2}{*}{\centering ML}
 & 16.45 & 10.31 & 5.22 & 2.82 & 41.87 & 26.85 & 14.11 & 7.57 \\ 
 &  &  & -1.52 & -1.27 & -0.96 & -0.46 & 0.98 & 0.99 & 0.99 & 1 \\ 
\cline{3-11} 
&&\multirow{2}{*}{\centering MLt}
 & 21.1 & 11.34 & 5.34 & 2.83 & 41.85 & 26.85 & 14.11 & 7.57 \\ 
 &  &  & 0.82 & -0.17 & -0.54 & -0.28 & 0.94 & 0.98 & 0.99 & 1 \\ 
\cline{3-11} 
&&\multirow{2}{*}{\centering MLt $10\sigma$}
 & 20.73 & 10.72 & 4.93 & 2.64 & 31.67 & 21.92 & 12.21 & 6.67 \\ 
 &  &  & 0.81 & -0.53 & -0.95 & -0.58 & 0.86 & 0.96 & 0.99 & 0.99 \\ 
\cline{3-11} 
&&\multirow{2}{*}{\centering MLt $100\sigma$}
 & 20.75 & 10.36 & 4.56 & 2.48 & 16.11 & 4.31 & 3.24 & 3.37 \\ 
 &  &  & 1.8 & 1.17 & 0.18 & -0.71 & 0.54 & 0.34 & 0.52 & 0.82 \\ 
\cline{3-11} 
\cline{2-11} 
\hline 
\multirow{12}{*}{\rotatebox[origin=c]{90}{\parbox[l]{3cm}{\centering $\eta = \sqrt{30}$ \  (low overlap) }}} 
&\multirow{6}{*}{\rotatebox[origin=c]{90}{\parbox[c]{3cm}{\centering $\sigma_{\varepsilon}=1.0$}}} 
&\multirow{2}{*}{\centering IPW}
 & 4.12 & 1.97 & 0.94 & 0.56 & 7.87 & 6.01 & 3.22 & 1.77 \\ 
 &  &  & -0.76 & -0.4 & -0.34 & -0.31 & 0.92 & 0.96 & 0.99 & 0.99 \\ 
\cline{3-11} 
&&\multirow{2}{*}{\centering AIPW}
 & 6.25 & 2.94 & 1.4 & 0.84 & 7.83 & 6 & 3.22 & 1.77 \\ 
 &  &  & -1.14 & -0.62 & -0.54 & -0.5 & 0.83 & 0.91 & 0.92 & 0.91 \\ 
\cline{3-11} 
&&\multirow{2}{*}{\centering OLS}
 & 4.04 & 2.81 & 2.53 & 2.46 & 6.53 & 3.79 & 1.84 & 0.93 \\ 
 &  &  & -2.36 & -2.43 & -2.43 & -2.44 & 0.81 & 0.74 & 0.19 & 0 \\ 
\cline{3-11} 
&&\multirow{2}{*}{\centering ML}
 & 2.91 & 2.01 & 1.27 & 0.81 & 4.96 & 3.19 & 2.09 & 1.45 \\ 
 &  &  & -2.18 & -1.78 & -1.18 & -0.77 & 0.77 & 0.84 & 0.94 & 0.99 \\ 
\cline{3-11} 
&&\multirow{2}{*}{\centering MLt}
 & 2.75 & 1.14 & 0.66 & 0.46 & 5.03 & 3.2 & 2.1 & 1.45 \\ 
 &  &  & 0.47 & -0.28 & -0.45 & -0.4 & 0.94 & 0.97 & 1 & 1 \\ 
\cline{3-11} 
&&\multirow{2}{*}{\centering MLt $10\sigma$}
 & 2.81 & 1.16 & 0.87 & 0.68 & 4.58 & 2.66 & 1.65 & 1.11 \\ 
 &  &  & 0.61 & -0.45 & -0.74 & -0.64 & 0.92 & 0.94 & 0.95 & 0.97 \\ 
\cline{3-11} 
&&\multirow{2}{*}{\centering MLt $100\sigma$}
 & 3.86 & 2.14 & 1.01 & 0.43 & 4.31 & 2.09 & 0.96 & 0.51 \\ 
 &  &  & 1.89 & 1.52 & 0.78 & -0.34 & 0.82 & 0.68 & 0.62 & 0.72 \\ 
\cline{3-11} 
\cline{2-11} 
&\multirow{6}{*}{\rotatebox[origin=c]{90}{\parbox[c]{3cm}{\centering $\sigma_{\varepsilon}=10.0$}}} 
&\multirow{2}{*}{\centering IPW}
 & 4.37 & 3.86 & 1.21 & 0.66 & 9.56 & 7.48 & 3.94 & 2.02 \\ 
 &  &  & -0.84 & -0.36 & -0.35 & -0.33 & 0.93 & 0.95 & 0.99 & 0.99 \\ 
\cline{3-11} 
&&\multirow{2}{*}{\centering AIPW}
 & 5.89 & 5.49 & 1.74 & 0.93 & 9.53 & 7.48 & 3.94 & 2.02 \\ 
 &  &  & -1.09 & -0.48 & -0.48 & -0.57 & 0.87 & 0.93 & 0.94 & 0.92 \\ 
\cline{3-11} 
&&\multirow{2}{*}{\centering OLS}
 & 4.78 & 3.19 & 2.66 & 2.54 & 9.03 & 5.15 & 2.38 & 1.2 \\ 
 &  &  & -2.3 & -2.49 & -2.5 & -2.5 & 0.93 & 0.87 & 0.46 & 0 \\ 
\cline{3-11} 
&&\multirow{2}{*}{\centering ML}
 & 3.3 & 2.29 & 1.46 & 0.86 & 6.2 & 4.21 & 2.62 & 1.73 \\ 
 &  &  & -2.28 & -1.8 & -1.25 & -0.76 & 0.89 & 0.91 & 0.93 & 0.98 \\ 
\cline{3-11} 
&&\multirow{2}{*}{\centering MLt}
 & 3.3 & 1.68 & 0.94 & 0.56 & 6.24 & 4.22 & 2.63 & 1.73 \\ 
 &  &  & 0.26 & -0.28 & -0.53 & -0.39 & 0.93 & 0.97 & 0.99 & 1 \\ 
\cline{3-11} 
&&\multirow{2}{*}{\centering MLt $10\sigma$}
 & 3.3 & 1.62 & 1.08 & 0.74 & 5.31 & 3.45 & 2.09 & 1.35 \\ 
 &  &  & 0.38 & -0.46 & -0.82 & -0.65 & 0.87 & 0.95 & 0.94 & 0.95 \\ 
\cline{3-11} 
&&\multirow{2}{*}{\centering MLt $100\sigma$}
 & 3.99 & 2.44 & 1.02 & 0.5 & 4.33 & 2.17 & 1 & 0.6 \\ 
 &  &  & 1.56 & 1.51 & 0.68 & -0.37 & 0.71 & 0.64 & 0.66 & 0.75 \\ 
\cline{3-11} 
\cline{2-11} 
&\multirow{6}{*}{\rotatebox[origin=c]{90}{\parbox[c]{3cm}{\centering $\sigma_{\varepsilon}=100.0$}}} 
&\multirow{2}{*}{\centering IPW}
 & 31.69 & 15.94 & 6.98 & 3.57 & 60.3 & 38.62 & 18.41 & 9.43 \\ 
 &  &  & -2.49 & -1.1 & -0.01 & -0.39 & 0.96 & 0.98 & 0.99 & 0.99 \\ 
\cline{3-11} 
&&\multirow{2}{*}{\centering AIPW}
 & 33.72 & 15.94 & 7.21 & 3.73 & 60.06 & 38.61 & 18.41 & 9.43 \\ 
 &  &  & -2.31 & -1.23 & -0.18 & -0.59 & 0.9 & 0.98 & 0.99 & 0.99 \\ 
\cline{3-11} 
&&\multirow{2}{*}{\centering OLS}
 & 29.68 & 13.09 & 6.03 & 3.74 & 67.56 & 34.27 & 15.36 & 7.7 \\ 
 &  &  & -2.38 & -3.07 & -2.14 & -2.51 & 0.97 & 0.99 & 0.99 & 0.96 \\ 
\cline{3-11} 
&&\multirow{2}{*}{\centering ML}
 & 16.52 & 9.96 & 5.77 & 3.5 & 40.93 & 27.46 & 15.84 & 9.4 \\ 
 &  &  & -2.46 & -2.05 & -1.01 & -0.81 & 0.98 & 0.99 & 0.99 & 0.99 \\ 
\cline{3-11} 
&&\multirow{2}{*}{\centering MLt}
 & 21.85 & 11.13 & 5.94 & 3.49 & 40.92 & 27.46 & 15.84 & 9.4 \\ 
 &  &  & -0.03 & -0.63 & -0.26 & -0.44 & 0.94 & 0.98 & 0.99 & 0.99 \\ 
\cline{3-11} 
&&\multirow{2}{*}{\centering MLt $10\sigma$}
 & 21.25 & 10.49 & 5.16 & 3.02 & 31.4 & 21.8 & 12.91 & 7.66 \\ 
 &  &  & 0.12 & -0.92 & -0.43 & -0.7 & 0.85 & 0.96 & 0.99 & 0.99 \\ 
\cline{3-11} 
&&\multirow{2}{*}{\centering MLt $100\sigma$}
 & 20.94 & 10.2 & 4.69 & 2.37 & 17.07 & 4.65 & 3.25 & 3.31 \\ 
 &  &  & 1.21 & 0.92 & 0.98 & -0.41 & 0.57 & 0.37 & 0.51 & 0.82 \\ 
\cline{3-11} 
\cline{2-11} 
\hline 
\end{tabular}
\\ 
}
\end{center}
\caption[Estimation error in the example of \citet*{hainmueller} with Outcome Design 3]{Root mean squared error (rmse), bias, and confidence interval half-width and coverage over 1000 replications
of Outcome Design 2 from \citet*{hainmueller}.
Here we take the tuning parameter $\sigma$ to be $0.1$ in the estimators ML and MLt. The notation MLt $10\sigma$ and $100\sigma$ indicates the substitution of $1$ and $10$ respectively.}
\end{figure}

\end{appendix}
\fi

\end{document}